\documentclass[a4paper]{article}
\usepackage{amsfonts}

\newtheorem{theorem}{Theorem}[section]
\newtheorem{lemma}[theorem]{Lemma}
\newtheorem{example}[theorem]{Example}
\newtheorem{proposition}[theorem]{Proposition}
\newtheorem{corollary}[theorem]{Corollary}
\newtheorem{definition}[theorem]{Definition}
\newtheorem{remark}[theorem]{Remark}
\newtheorem{lemexp}[theorem]{Lemma and Example}
\newcommand{\id}{{\rm id}}
\newcommand{\CB}{{\rm CB}}

\begin{document}

\title{Cohomology of Hopf $C^{\ast }$-algebras and Hopf von Neumann algebras}
\author{Chi-Keung Ng}
\date{}
\maketitle

\begin{abstract}
We will define two canonical cohomology theories for Hopf
$C^{\ast }$-algebras and for Hopf von Neumann algebras (with coefficients in their
comodules). We will then study the situations when these cohomologies
vanish. The cases of locally compact groups and compact quantum groups will
be considered in more details.
\end{abstract}

\noindent {\small 1991 AMS Mathematics Classification number: Primary:
46L55, 46L05; Secondary: 43A07, 22D25}

\bigskip

\bigskip

\bigskip

\bigskip

\bigskip

\bigskip

\bigskip

\bigskip

\bigskip

\bigskip

\emph{The first statement of Remark 1.9(a) in the original published article does not follow directly from 
Lemma 1.7(a) but it does follow very easily from the argument of Proposition 1.10. 
Therefore, we change the presentations of Remark 1.9 and Proposition 1.10 (one line was
added in the proof of Proposition 1.10). 
Please find in the following the corrected version of this paper
(any further question, comment or correction is welcomed).}

\pagebreak

\section*{0\hspace{0.25in}Introduction}

\bigskip

Cohomology theory is an important subject in many branches of Mathematics.
In the field of Banach algebras, the vanishing of cohomology defines the
interesting notion of amenability which, in the case of group algebras of
locally compact groups, generalised the concept of amenable groups. In \cite
{Ruan2}, Ruan studied the operator cohomology of completely contractive
Banach algebras and defined the notion of operator amenability. He showed
that the Fourier algebra of a locally compact group is operator amenable if
and only if the group is amenable. The objective of this paper is to define
and study some cohomology theories of Hopf $C^{\ast }$-algebras as well as
Hopf von Neumann algebras.
(Note that the notion
of Hopf $C^*$-algebras that we use here is the same as that in \cite[0.1]{BS}
even though for most of the cases, we will also assume the extra condition
mentioned in \cite[0.1]{BS} -- which is called {\em saturated} in this paper;
see Definition \ref{Hopf}).
These will be important tools for the study of
``locally compact quantum groups''. In particular, we will also investigate the
situation when these cohomologies vanish.

\medskip

In fact, there are four natural approaches to define cohomology theory for
Hopf $C^{\ast }$-algebras and Hopf von Neumann algebras:

\begin{enumerate}
\item  analogue of the deformation cohomology for Hopf algebras (see e.g. \cite
{Shn});

\item  analogue of the cyclic cohomology for Hopf algebras (see e.g. \cite{CM});

\item  generalisation of the group cohomology (see e.g. \cite{Pier}, \cite{Evens}
and \cite{Wilson} for three different meanings of this);

\item  ``dual analogue'' of the Banach algebra homology (see. e.g. \cite{John}).
\end{enumerate}
\noindent
For the first approach, we note that the definition of the cochain complex of
the deformation cohomology for Hopf algebras requires operations involving
both the product and the coproduct in a way which is unnatural for operator
algebras (in particular, the maps are not bounded under the ``default tensor
product''). However, by using some technical results concerning the extended
Haagerup tensor product (see \cite{EKR}), we can still define this sort of
analogue for Hopf von Neumann algebras. The same difficulty (together with
some other) arise in the case of cyclic cohomology.
(We didn't try to find
this analogue but even if this can be done, it is believed that the analogue
can only be defined for Hopf von Neumann algebras -- with the help of the
reshuffle map in \cite{EKR}).
Nevertheless, we will not study
these in this paper. Instead, we will study the kind of cohomology theories
which are related to those of the third and the fourth approaches.

\medskip

More precisely, the \emph{natural cohomology} in this paper is defined
according to the fourth approach (a comparison of this cohomology theory with
the existing theories of ``group cohomology'' can be found in Example
\ref{B17}) whereas the \emph{dual cohomology} is conceptually a dual
version of the natural cohomology (although not directly related). Note
that, even in the case of locally compact groups, the dual cohomology
is different from any kind of the cohomology theories known so far (a
comparison of the first dual cohomology theory with other cohomology theories
can be found in Remark \ref{4.7}(a) and a comparison of the dual
cohomology theory with a cohomology theory of coalgebras can be found in
Remark \ref{rmstrcd}(c)). Moreover, we know that for a locally compact
group $G$, the cohomology theory for the Fourier algebra $A(G)$ as a
Banach algebra is different from the cohomology theory for $A(G)$ as a
completely contractive Banach algebra (see \cite{John2} and
\cite{Ruan2}) even though they are formally defined in the same way.
Hence there is no reason to believe that the dual cohomology of
$C_{0}(G)$ or $C_{r}^{\ast }(G)$ will behave as either the Banach
algebra cohomology or the operator cohomology of $A(G)$ (or
$L^{1}(G)$). Therefore, these cohomology theories deserve detail
study.

\medskip

Before we can define the cohomology theories, we
need to understand first of all, the comodules of Hopf $C^{\ast }$-algebras
(and Hopf von Neumann algebras). A comodule of a Hopf $C^{\ast }$-algebra $S$
can be thought of as a vector space $X$ together with some topological
structures as
well as a ``scalar-comultiplication'' $\beta $ from $X$ to a kind of
topological completion of the algebraic tensor product $X\odot S$. Moreover,
this should include the case of coactions on $C^{\ast }$-algebras. This
means that the range of the coaction should be the ``space of multipliers''
of the completion of $X\odot S$. In fact, if we want to define a comodule
structure on a Banach space $X$, we will first come across the problem of
getting a right topology for $X\odot S$. Furthermore, the set of
``multipliers'' does not behave nicely as required. There are also some
other technical difficulties. However, if we consider operator spaces
instead of Banach spaces, all these difficulties can be overcome.

\medskip

In section 1 of this paper, we will recall some basic materials about
operator spaces. We will then study multipliers of operator bimodules and
give the definition as well as examples of comodules of Hopf
$C^{\ast }$-algebras.
Note that in most of the cases in this paper (except for some
results in the final section as well as in the appendix), we will assume
that the Hopf $C^*$-algebras are saturated (see Definition \ref{Hopf}).

\medskip

In the second section, we will define the natural and the dual cohomology
theories mentioned above. We will study some situations when these cohomology
theories vanish. If the Hopf $C^{\ast }$-algebra is saturated and unital,
(i.e. a compact quantum group; see
\cite{VW} or \cite{Wor}), we will give some interesting
equivalent conditions for the vanishing of the (two-sided) dual cohomology
(Theorem \ref{B20}). We will then show (in Corollary \ref{discrete}) that in
the case when the Hopf $C^*$-algebra is $C_{r}^{\ast }(\Gamma )$
where $\Gamma $ is a discrete group, the
vanishing of the dual cohomology is equivalent to the amenability of
$\Gamma$.
More generally, we will see in the final section that for a general
locally compact group $G$, the vanishing of the first dual cohomology of
$C_{r}^{\ast }(G)$ is equivalent to the amenability of $G$
(Theorem \ref{D9}(b)).
On the other hand, the vanishing of the first dual cohomology of
$C_{0}(G)$ is again equivalent to $G$ being amenable (Theorem \ref{D9}(a)).
These are surprising since the dual cohomology theory of Hopf
$C^{\ast }$-algebras are thought to be not as sharp as theirs von Neumann algebra
counterpart (see Remark \ref{4.7}(a)).

\medskip

In section three, we will define comodules of Hopf von Neumann algebras and
give the analogues of the above cohomology theories in this case.
In particular, we will
show that there is a natural Hopf von Neumann algebra comodule structure on
the dual space of a comodule and the dual cohomology with coefficients in
a given comodule is the same as the natural cohomology with coefficient in
the corresponding dual comodule (which is not true in the Hopf
$C^{\ast }$-algebra case). Furthermore, the dual cohomology also coincides with the
operator cohomology (see \cite{Ruan2}) of the predual of the underlying Hopf
von Neumann algebra.

\medskip

In the final section, we will give some interesting relations between the
vanishing of the (one sided and two sided) dual cohomologies and
amenability. In particular, we will give a characterisation of the
amenability of discrete semi-groups in terms of the dual cohomology theory
and will prove the characterisation of amenable locally compact groups
mentioned above. We will also consider the amenability of general Hopf $C^{\ast
}$-algebras (see \cite{Ng2}).

\bigskip

\noindent Acknowledgement

\medskip

\emph{The author would like to thank Prof. Z. J. Ruan for reading a very
early version of this article and for the comments. The author would
also like to thank those talked to him about this subject during
the EU conference on }$C^{\ast }$-algebras and \emph{Non-commutative
Geometry in Copenhagen 1998, especially those suggested him to
look at the two-sided case. }

\bigskip

\bigskip

\section{Operator modules and coactions on Operator spaces}

\bigskip

In this section, we will recall some properties of operator spaces and
define coactions on them.

\bigskip
\noindent
{\bf Notation:}
{\em Throughout this paper, $X$, $Y$ and $Z$ are operator
spaces. $\odot $ is the algebraic tensor product while $\otimes $ is the
operator spatial tensor product and $\hat{\otimes}$ is the operator
projective tensor product (see e.g. \cite{BP}).}

\medskip
\begin{definition}
\label{1.e} Let $B$ be a normed $\ast $-algebra and $N$ be a normed space.

\smallskip
\noindent (a) A norm $\Vert \cdot \Vert _{\alpha }$ on the algebraic tensor
product $B\odot N$ is said to be a {\em $B$-bimodule cross norm} if it is a
cross norm (i.e. $\|a\otimes x\| = \|a\|\|x\|$) and $\Vert a\cdot
z\cdot b\Vert _{\alpha }\leq \Vert a\Vert _{B}\Vert z\Vert _{\alpha}
\Vert b\Vert_{B}$ for any $z\in B\odot N$ and $a,b\in B$
(where $a\cdot (c\otimes n)\cdot b=acb\otimes n$).

\smallskip
\noindent (b) A $B$-bimodule cross norm on $B\odot N$ is said to be a 
{\em $L^{\infty }$ $B$-bimodule cross norm} if for any
disjoint (self-adjoint) projections
$p,q\in B$ (i.e. $pq=0$), $\Vert p\cdot z\cdot p+q\cdot z\cdot q\Vert
_{\alpha }=\max \{\Vert p\cdot z\cdot p\Vert _{\alpha },\Vert q\cdot
z\cdot q\Vert _{\alpha }\}$.
\end{definition}

\medskip
The following proposition is probably well known and a proof of which can be
found in \cite{Ng3} (see also \cite{Muhly}).

\medskip
\begin{proposition}
\label{1.f} (a) For any Banach space $N$, there is an one to one
correspondence between the operator space structures on $N$ and the
$L^{\infty }$ $\cal K$-bimodule cross norms on $\mathcal{K}\odot N$ (where $\cal K$
is the space of compact operators on the separable Hilbert space $l^2$).
In this case, the $L^{\infty }$ $\cal K$-bimodule cross norm on $\mathcal{K}\odot
N$ is given by the
operator spatial tensor product $\mathcal{K}\otimes N$.

\smallskip
\noindent (b) A linear map $S$ from $X$ to $Y$ is completely bounded if and
only if the map $\mathrm{id}_{\mathcal{K}}\otimes S$ extends to a bounded
map from $\mathcal{K}\otimes X$ to $\mathcal{K}\otimes Y$. In this case,
$\Vert S\Vert _{\mathrm{cb}}=\Vert \mathrm{id}_{\mathcal{K}}\otimes S\Vert $.

\smallskip
\noindent (c) $\mathcal{L_{K}}(\mathcal{K}\otimes X;\mathcal{K}\otimes Y)=
\{\mathrm{id}_{\mathcal{K}}\otimes S:S\in \mathrm{CB}(X;Y)\}$ (where
$\mathcal{L_{K}}$ means the set of all bounded $\mathcal{K}$-bimodule maps).
Consequently, $\mathrm{CB}(X;Y)\cong \mathcal{L_{K}}(\mathcal{K}\otimes X;
\mathcal{K}\otimes Y)$ as normed spaces.

\smallskip
\noindent (d) The canonical injection from $\mathcal{K}\odot \mathrm{CB}(X;Y)
$ to $\mathrm{CB}(X;\mathcal{K}\otimes Y)$ gives the natural operator space
structure on $\mathrm{CB}(X;Y)$.
\end{proposition}

\medskip
The following are some easy facts about projective tensor product.
Note that part (a) is easy to check while parts (b) and (c) are the ideas
behide the
definition of operator projective tensor product (see \cite[5.4]{BP}
as well as the paragraph following \cite[5.3]{BP}).
Moreover, part (d) is a direct consequence of part (c).

\medskip
\begin{lemma}
\label{A25} Let $X$, $Y$ and $Z$ be operator spaces.

\smallskip
\noindent (a) For any $T\in \mathrm{CB}(X;Z)$, the map $T\otimes \mathrm{id}$
on the algebraic tensor product extends to a completely bounded map from
$X\hat{\otimes}Y$ to $Z\hat{\otimes} Y$.

\smallskip
\noindent (b) $\mathrm{CB}(X\hat{\otimes}Y;Z)\cong \mathrm{CB}(X;\mathrm{CB}
(Y;Z))\cong \mathrm{JCB}(X\times Y;Z)$ (jointly completely bounded bilinear
maps from $X\times Y$ to $Z$; see \cite{BP}).

\smallskip
\noindent (c) $\mathrm{CB}(X;Y^{\ast })=(X\hat{\otimes}Y)^{\ast }$.

\smallskip
\noindent (d) $\mathrm{CB}(X\hat{\otimes}Y;Z^{\ast })\cong \mathrm{CB}(Y
\hat{\otimes}Z;X^{\ast })$.
\end{lemma}

\medskip
\begin{remark}
\label{A26} We call the identification in Lemma \ref{A25}(b) (i.e.
$\mathrm{CB}(X\hat{\otimes}Y;Z)\cong \mathrm{CB}(X;\mathrm{CB}(Y;Z))$)
the \emph{standard identification} whereas the identification
$\mathrm{CB}(Y\hat{\otimes}X;Z)\cong \mathrm{CB}(X;\mathrm{CB}(Y;Z))$
will be called the
\emph{reverse identification}. This distinction is important
when $X=Y$ (in which case the two identifications look the same but have
different meanings; see the paragraph after Definition \ref{C5}).
\end{remark}

\medskip
The following can be found in \cite[1.5(a)]{Ng3} and is again well known.

\medskip
\begin{lemma}
\label{A1} Let $X$ be a closed subspace of an operator space $Y$.
If $\psi $ is a completely
bounded map from $Y$ to $Z$ and $X\subseteq \mathrm{Ker}(\psi )$, then $\psi
$ induces a completely bounded map $\hat{\psi}$ from $Y/X$ to $Z$ such that
$\Vert \hat{\psi}\Vert _{cb}\leq \left\| \psi \right\| _{cb}$ and $\hat{\psi}
\circ q =\psi$ (where $q$ is the canonical map from $Y$ to $Y/X$).
\end{lemma}

\medskip
The following trivial lemma sets up some notations to be
used later.

\medskip
\begin{lemma}
\label{A27} Let $W$, $X$, $Y$ and $Z$ be operator spaces.

\smallskip
\noindent (a) Any element $F\in \mathrm{CB}(Z;W)$ induces a completely
bounded map $\tilde{F}$ from $\mathrm{CB}(Y;Z)$ to $\mathrm{CB}(Y;W)$ such
that $\tilde{F}(T)=F\circ T$ and $\Vert \tilde{F}\Vert _{\mathrm{cb}}\leq
\Vert F\Vert _{\mathrm{cb}}$.

\smallskip
\noindent (b) For any $T\in \mathrm{CB}(X;\mathrm{CB}(Y;Z))$, there is a
completely bounded map $T^{\#}:\mathrm{CB}(Z;W)\longrightarrow \mathrm{CB}(X;
\mathrm{CB}(Y;W))$ (respectively, $T^{0}:\mathrm{CB}(W;Y)\longrightarrow
\mathrm{CB}(X;\mathrm{CB}(W;Z))$) such that $T^{\#}(F)(x)(y)=F(T(x)(y))$ and
$\left\| T^{\#}\right\| _{\mathrm{cb}}\leq \left\| T\right\| _{\mathrm{cb}}$
(respectively, $T^{0}(F)(x)(w)=T(x)(F(w))$ and
$\left\| T^{0}\right\| _{\mathrm{cb}}\leq \left\| T\right\| _{\mathrm{cb}}$).
\end{lemma}

\medskip
We would like to study multipliers of operator $A$-bimodules for
a $C^*$-algebra $A$.
Let us first look at the multipliers of Banach $A$-bimodules.
Given an $A$-bimodule $N$, let $M_{A}^{l}(N)$ (respectively, $M_{A}^{r}(N)$)
be the set of all linear maps from $A$ to $N$ that respect the right
(respectively, left) $A$-multiplications, i.e. the set of all left
(respectively, right) multipliers. Let $M_{A}(N)=\{(l,r)\in
M_{A}^{l}(N)\times M_{A}^{r}(N):a\cdot l(b)=r(a)\cdot b$ \textrm{for
any }$a,b\in A\}$. Elements in $M_{A}(N)$ are called the
\emph{multipliers} of $N$. A bimodule $N$ is said to be
\emph{essential} if both $A\cdot N$ and $N\cdot A$ are dense in $N$.
If $A$ is unital, $N$ is essential simply means that $N$ is a unital 
$A$-bimodule.
Moreover, in this case, $M_A(N)=M^l_A(N)=M^r_A(N)=N$.

\bigskip
\noindent
{\bf Notation:}
{\em From now on, until the end of this section, $A$ is a $C^{\ast }$-algebra.
Moreover, if $(l,r)\in M_A(N)$ and $a\in A$, we will denote $a\cdot (l,r) = r(a)$
and $(l,r)\cdot a = l(a)$.}

\medskip
\begin{lemma}
\label{1.a} Let $A$ be a $C^{\ast }$-algebra and $N$ be an essential $A$-bimodule.

\smallskip
\noindent (a) Any left or right multiplier on $N$ is automatically bounded.

\smallskip
\noindent (b) For any $(l,r)\in M_{A}(N)$, we have $\Vert l\Vert
_{M_{A}^{l}(N)}=\Vert r\Vert _{M_{A}^{r}(N)}$.

\smallskip
\noindent (c) $M_{A}(N)$ is a Banach space for the norm defined by
$\Vert (l,r) \Vert = \Vert l\Vert
_{M_{A}^{l}(N)}=\Vert r\Vert _{M_{A}^{r}(N)}$.
\end{lemma}

\medskip
In fact, part (a) follows from a similar argument as in \cite[3.12.2]{Ped}
and part (b) follows from the fact that $A$ has an approximate unit for $N$
while part (c) is more or less obvious.

\medskip

Suppose that $Y$ is an operator space with an $A$-bimodule structure.
Then $Y$ is called an \emph{operator }$A$\emph{-bimodule} if $
\mathcal{K}\otimes Y$ is a $\mathcal{K}\otimes A$-bimodule. In this case, if
$Y$ is essential as an $A$-bimodule, then $\mathcal{K}\otimes Y$ is an
essential $\mathcal{K}\otimes A$-bimodule and we call $Y$ an
\emph{essential operator }$A$\emph{-bimodule}. Moreover, by using
\cite[3.3]{CES} and some simple arguments concerning the essentialness
of the bimodule $Y$ as well as employing the trick of replacing $Y$ with
$\mathcal{K\otimes }Y$, we have the following representation lemma
(a detail argument can be found in \cite{Ng3}).
The triple $(\phi ,\pi ,\psi )$ satisfying the relation in this lemma
is called a \emph{spatial realisation }of $Y$.

\medskip
\begin{lemma}
\label{1.b} Let $Y$ be an essential operator $A$-bimodule. Then there exist
Hilbert spaces $H$ and $K$ as well as a complete isometry $\pi $ from $Y$ to $
\mathcal{L}(H;K)$ and faithful non-degenerate representations $\psi $ and $
\phi $ of $A$ on $\mathcal{L}(H)$ and $\mathcal{L}(K)$ respectively such
that $\phi (b)\pi (y)\psi (a)=\pi (b\cdot y\cdot a)$ for all $a,b\in A$ and $
y\in Y$.
\end{lemma}

\medskip
This suggests another way to define multipliers: $M_{A}^{\pi }(Y)=\{m\in
\mathcal{L}(H;K):\phi (A)m,m\psi (A)\subseteq \pi (Y)\}$. However, it is not
obvious that this definition is independent of the choice of the spatial
realisation. Nevertheless, we will see later that it is indeed completely
isometrically isomorphic to $M_{A}(Y)$ (regarded as an operator $A$-bimodule). Let us first give a natural operator space structure on $
M_{A}(Y) $.

\medskip
\begin{remark}
\label{1.c} 
(a) Suppose that $M_{A,cb}^l(Y) = M_A^l(Y)\cap \CB(A;Y)$ and $M_{A,cb}^r(Y) = M_A^r(Y)\cap \CB(A;Y)$
with the induced operator space structures. 
More precisely, the operator space structure on $M_{A,cb}^l(Y)$
is given by the canonical injection from $\mathcal{K}\odot
M_{A,cb}^{l}(Y)$ to $M_{A,cb}^{l}(\mathcal{K}\otimes Y)\subseteq \mathrm{CB}(A;
\mathcal{K}\otimes Y)$ (see Proposition \ref{1.f}(d)). 
We denote by $\Vert \cdot \Vert _{\mathrm{usu}}$ and 
$\Vert \cdot \Vert $ the norms on $M_{A}^{l}(Y)$ and $M_{A,cb}^{l}(Y)$ induced from 
$\mathcal{L}(A;Y)$ and $\CB(A;Y)$ respectively.
Since $\mathrm{CB}(A;Y)$ can be regarded as the subspace $\mathcal{L_{K}}(
\mathcal{K}\otimes A;\mathcal{K}\otimes Y)$ of $\mathcal{L}
(\mathcal{K}\otimes A;\mathcal{K}\otimes Y)$
(see Proposition \ref{1.f}(b) and (c)), the
canonical map from $(M_{A,cb}^{l}(Y),\Vert \cdot \Vert )$ to $(M_{\mathcal{K}
\otimes A}^{l}(\mathcal{K}\otimes Y),\Vert \cdot \Vert _{\mathrm{usu}})$ is
an isometry. 
Therefore, the operator space structure on $M_{A,cb}^{l}(Y)$ is
given by the canonical embedding from $\mathcal{K}\odot M_{A,cb}^{l}(Y)$ to 
$(M_{\mathcal{K}\otimes A}^{l}(\mathcal{K}\otimes \mathcal{K}\otimes Y),\Vert
\cdot \Vert _{\mathrm{usu}})$ (as $(M_{A,cb}^{l}(\mathcal{K}\otimes Y), \Vert \cdot \Vert )$ 
can be regarded as its subspace).
The same is true for $M_{A,cb}^r(Y)$.

\smallskip
\noindent 
(b) 
Lemma \ref{1.a}(b) implies that $(M_{\mathcal{K}\otimes A}(
\mathcal{K}\otimes Y),\Vert \cdot \Vert _{\mathrm{usu}})$ is simultaneously a norm subspace of both 
$(M_{\mathcal{K}\otimes A}^{l}(\mathcal{K}\otimes Y),\Vert \cdot \Vert _{\mathrm{usu}})$ and 
$(M_{\mathcal{K}\otimes A}^{r}(\mathcal{K}\otimes Y),\Vert \cdot \Vert _{\mathrm{usu}})$. 
Thus, the norms induced on $M_{A,cb}(Y) = M_A(Y)\cap (\CB(A;Y)\times \CB(A;Y))$ from 
$(M_{A}^{l}(Y),\Vert \cdot \Vert )$ and $(M_{A}^{r}(Y),\Vert \cdot \Vert )$ coincide. 
Similarly, $M_{\mathcal{K}\otimes A}(\mathcal{K}\otimes \mathcal{K}\otimes Y)$ 
is simultaneously a norm subspace of both 
$(M_{\mathcal{K}\otimes A}^{l}(\mathcal{K}\otimes \mathcal{K}\otimes Y),\Vert \cdot \Vert _{\mathrm{usu}})$ and 
$(M_{\mathcal{K}\otimes A}^{r}(\mathcal{K}\otimes \mathcal{K}\otimes Y),\Vert \cdot \Vert _{\mathrm{usu}})$. 
Therefore, the two embeddings from $\mathcal{K}\odot M_{A,cb}(Y)$ to $\mathcal{L}_{\cal K}(\mathcal{K}\otimes A; 
\mathcal{K}\otimes \mathcal{K}\otimes Y)$ (induced from $M^l$ and $M^r$)
give the same operator space structure on $M_{A,cb}(Y)$ and 
we use this structure by default. 
\end{remark}

\medskip

\begin{proposition}
\label{1.d} 
Let $A$ be a $C^*$-algebra and $Y$ be an essential
operator $A$-bimodule. For any spatial realisation ($\phi ,\pi ,\psi )$ of $Y$
(see the statement before Lemma 1.8), there exists an isometry
$\Psi $ from $(M_{A}(Y), \Vert \cdot \Vert _{\mathrm{usu}})$ onto $M_{A}^{\pi }(Y)$ such that 
$\Psi(l,r)\psi (a)=\pi (l(a))$ and $\phi (a)\Psi (l,r)=\pi (r(a))$ 
($a\in A$; $(l,r)\in M_{A}(Y)$). 
Moreover, $(M_{A,cb}(Y), \Vert \cdot \Vert) = (M_A(Y), \Vert \cdot \Vert _{\mathrm{usu}})$ and 
$\Psi$ is a complete isometry. 
\end{proposition}
\noindent \textbf{Proof:} Let $\{a_{i}\}$ be an approximate unit of $A$.
Suppose that $(l,r)\in M_{A}(Y)$.
The net $\{\pi (l(a_{i}))\}$ converges strongly to
an element $m\in \mathcal{L}(H;K)$ (as $\psi $ is non-degenerate and $l$ is
bounded). It is clear that $m\psi (a)=\pi (l(a))$ and $\phi (b)m\psi (a)=\pi
(b\cdot l(a))=\pi (r(b))\psi (a)$ for any $a,b\in A$. Moreover,
\[
\Vert m\Vert =\sup \{\Vert m\psi (a)\Vert :a\in A;\Vert a\Vert \leq
1\}=\Vert l\Vert _{\mathrm{usu}}=\Vert (l,r)\Vert _{\mathrm{usu}}.
\]
Hence it is not hard to see that the map $\Psi $ that sends $(l,r)$ to $m$
is a surjective isometry from $(M_{A}(Y),\Vert \cdot \Vert _{\mathrm{usu}})$
to $M_{A}^{\pi }(Y)$ (note that any element in $M_{A}^{\pi }(Y)$ defines in
the obvious way, an element in $M_{A}(Y)$) which satisfies the required
equalities. 
As a consequence, $M_{A,cb}(Y) = M_A(Y)$ (because both $a\mapsto \pi^{-1}(m\psi(a))$ 
and $a\mapsto \pi^{-1}(\phi(a)m)$ are 
completely bounded maps for $m\in M_{A}^{\pi }(Y)$). 
It remains to show that $\Vert \cdot \Vert _{\mathrm{usu}}$
coincides with $\Vert \cdot \Vert $ and $\Psi $ is a complete isometry.
Observe that by replacing $Y$ with $\mathcal{K}\otimes Y$ and $A$ with $
\mathcal{K}\otimes A$, we have an isometry $\Psi ^{\prime }$ from $(M_{
\mathcal{K}\otimes A}(\mathcal{K}\otimes Y),\Vert \cdot \Vert _{\mathrm{usu}
})$ to $M_{\mathcal{K}\otimes A}^{\mathrm{id}\otimes \pi }(\mathcal{K}
\otimes Y)$. For any $k\in \mathcal{K}$ and $a\in A$,
\[
\Psi ^{\prime }(\mathrm{id}_{\mathcal{K}}\otimes l,\mathrm{id}_{\mathcal{K}
}\otimes r)(k\otimes a)=(1\otimes \Psi (l,r))(k\otimes a).
\]
Thus, $\Psi $ is an isometry from $(M_{A}(Y),\Vert \cdot \Vert )$ to $
M_{A}^{\pi }(Y)$ (recall from Remark \ref{1.c} that $\Vert (l,r)\Vert =\Vert
(\mathrm{id}_{\mathcal{K}}\otimes l,\mathrm{id}_{\mathcal{K}}\otimes r)\Vert
_{\mathrm{usu}}$). This also shows that $\Vert \cdot \Vert $ and $\Vert
\cdot \Vert _{\mathrm{usu}}$ agree on $M_{A}(Y)$. Now, if we replace $Y$ with $
\mathcal{K}\otimes Y$ only, we have an isometry $\Phi $ from $(M_{A}(
\mathcal{K}\otimes Y),\Vert \cdot \Vert )$ to $M_{A}^{\mathrm{id}\otimes \pi
}(\mathcal{K}\otimes Y)$. For any $k,k^{\prime }\in \mathcal{K}$ and $a\in A$
, $k\otimes (l,r)$ can be regarded as an element of $M_{A}(\mathcal{K}
\otimes Y)$ and
\[
\Phi (k\otimes (l,r))(k^{\prime }\otimes a)=(k\otimes \Psi (l,r))(k^{\prime
}\otimes a).
\]
Therefore, the map $\mathrm{id}_{\mathcal{K}}\otimes \Psi $ from $\mathcal{K}\otimes
M_{A}(Y)$ to $\mathcal{K}\otimes M_{A}^{\pi }(Y)\subseteq M_{A}^{\mathrm{id}
\otimes \pi }(\mathcal{K}\otimes Y)$ is an isometry (note that $(\mathcal{K}
\otimes M_A(Y), \Vert\cdot\Vert)$ is a subspace of $(M_A(\mathcal{K}
\otimes Y), \Vert\cdot\Vert)$ by Remark \ref{1.c})
and hence $\Psi $ is a complete isometry by Proposition \ref{1.f}(b).

\medskip
\begin{corollary}
\label{1.11} (a) $Y$ is an operator subspace of $M_{A}(Y)$.

\smallskip
\noindent (b) $M_{A}(Y)$ is a unital operator $M(A)$-bimodule.

\smallskip
\noindent (c) If $B$ is another $C^{\ast }$-algebra and $Z$ is an essential
operator $B$-bimodule, then there exists a complete isometry from $
M_{B}(M_{A}(Y)\otimes Z)$ to $M_{A\otimes B}(Y\otimes Z)$ that respects both
the $A$-bimodule and the $B$-bimodule structures.
\end{corollary}

\medskip
\noindent
{\bf Notation:}
{\em From now on, we may use the identification in Proposition \ref{1.d}
implicitly and will regard $M_{B}(M_{A}(Y)\otimes Z)$ as subspace of $
M_{A\otimes B}(Y\otimes Z)$.}

\medskip
\begin{lemma}
\label{1.12} (a) Let $X$ and $Y$ be essential operator $A$-bimodules.
Suppose that $\varphi $ is a completely bounded $A$-bimodule map from $X$ to
$Y$. Then $\varphi $ induces a completely bounded $M(A)$-bimodule map, again
denoted by $\varphi $, from $M_{A}(X)$ to $M_{A}(Y)$. If $\varphi $ is
completely isometric, then so is the induced map.

\smallskip
\noindent (b) Let $A$ and $B$ be $C^{\ast }$-algebras and $\psi
:A\longrightarrow M(B)$ be a non-degenerate $\ast $-homomorphism. Then $
\mathrm{id}_{Z}\otimes \psi $ extends to a complete contraction $\mathrm{id}
_{Z}\otimes \psi :M_{A}(Z\otimes A)\longrightarrow M_{B}(Z\otimes B)$ such
that $(\mathrm{id}_{Z}\otimes \psi )(m\cdot a)=(\mathrm{id}_{Z}\otimes \psi
)(m)\cdot \psi (a)$ and $(\mathrm{id}_{Z}\otimes \psi )(a\cdot m)=\psi
(a)\cdot (\mathrm{id}_{Z}\otimes \psi )(m)$ ($m\in M_{A}(Z\otimes A)$; $a\in
M(A)$). If $\psi $ is injective, then $\mathrm{id}_{Z}\otimes \psi $ is a
complete isometry. Furthermore, if $\phi $ is a completely bounded map from $
Z$ to another operator space $Z^{\prime }$, then $(\phi \otimes \mathrm{id})(
\mathrm{id}\otimes \psi )=(\mathrm{id}\otimes \psi )(\phi \otimes \mathrm{id}
)$ on $M_{A}(Z\otimes A)$.

\smallskip
\noindent (c) For any $g\in A^{\ast }$ and $T\in \mathrm{CB}(X;Y)$
($X$ and $Y$ are operator spaces), we
have $T\circ (\mathrm{id}\otimes g)=(\mathrm{id}\otimes g)(T\otimes \mathrm{id})$
on $M_{A}(X\otimes A)$.
\end{lemma}

\medskip
The map in part (a) is induced by the completely bounded map given in
Lemma \ref{A27}(a).
The first two statements of part (b) follow from
Proposition \ref{1.d} while the last statement follows from the fact that $
(\phi \otimes \mathrm{id})(m)(1\otimes b)=(\phi \otimes \mathrm{id}
)(m(1\otimes b))$ under the identification in Proposition \ref{1.d}.
The map $\id\otimes g$ in part (c) is defined on $M_A(X\otimes A)
\subseteq M(\mathcal{K}(H)\otimes A)$ when $X\subseteq \mathcal{K}(H)$.
It is well defined and satisfies the equality in (c) because
$g$ can be decomposed as $a\cdot g^{\prime }$ where $a\in A$ and
$g^{\prime }\in A^{\ast }$.

\medskip

Next, we will recall from \cite[0.1 \& 0.2]{BS} the notion of
Hopf $C^{\ast }$-algebras and their coactions
(even though we change some of the terminology in our translation).

\medskip
\begin{definition}
\label{Hopf}
(a) Let $S$ be a $C^{\ast }$-algebra with a non-degenerate
$\ast $-homomorphism $\delta $ from $S$ to $M(S\otimes S)$.
Then $(S,\delta)$ is said to be a {\em Hopf $C^{\ast }$-algebra} if $\delta (S)(1\otimes S),\delta
(S)(S\otimes 1)\subseteq S\otimes S$ and $(\delta \otimes \mathrm{id})\delta
=(\mathrm{id}\otimes \delta )\delta $. In this case, $\delta $ is called a
{\em coproduct} of $S$. Moreover, a Hopf $C^{\ast }$-algebra
$(S,\delta)$ is said to be {\em saturated} if both of the vector spaces
$\delta (S)(1\otimes S)$ and
$\delta (S)(S\otimes 1)$ are dense in $S\otimes S$.

\smallskip
\noindent (b) Let $\Re$ be a von Neumann algebra with a
weak*-continuous unital $\ast $-homomorphism from $\Re$ to
$\Re\bar{\otimes}\Re$.
Then $(\Re,\delta)$ is said to be a {\em Hopf von
Neumann algebra} if $(\delta \otimes \mathrm{id})\delta =(\mathrm{id}\otimes \delta
)\delta $. A Hopf von Neumann algebra $\Re$ is said to be
{\em saturated} if both $\delta (\Re)(1\otimes \Re)$ and $\delta (\Re
)(\Re\otimes 1)$ are weak*-dense in $\Re\bar{\otimes}
\Re$.

\smallskip
\noindent
(c) Let $A$ be a $C^*$-algebra and $\cal M$ be a von Neumann algebra.
A non-degenerate $*$-homomorphism $\beta$ from $A$ to $M_S(A\otimes S)$
(which is a $C^*$-algebra) is said be a {\em coaction of $S$ on $A$}
if $(\beta\otimes \id)\circ \beta = (\id\otimes
\delta)\circ \beta$.
Similarly, a normal $*$-homomorphism $\beta$ from $\cal M$ to
$\mathcal{M}\bar\otimes \Re$ is said to be {\em a coaction of $\Re$
on $\cal M$} if $(\beta\otimes \id)\circ \beta = (\id\otimes
\delta)\circ \beta$.
\end{definition}

\medskip
\noindent
{\bf Notation:}
{\em From now on, unless specified, $(S,\delta)$ is a saturated
Hopf $C^{\ast }$-algebra and $(\Re,\delta)$ is a saturated Hopf von
Neumann algebra.
Moreover, when there is no confusion arise, we will simply use $S$ and
$\Re$ to denote the Hopf $C^*$-algebra and Hopf von Neumann algebra respectively.}

\medskip
\begin{definition}
\label{2.a}
Suppose that $(R,\delta)$ is a (not necessarily saturated) Hopf $C^*$-algebra.

\smallskip\noindent
(a) Let $\beta $ be a completely bounded map from $X$ to
$M_{R}(X\otimes R)$. Then $\beta $ is said to be a \emph{right
coaction} of
$R$ on $X$ if $(\beta \otimes \mathrm{id})\beta =(\mathrm{id}
\otimes \delta )\beta \in \mathrm{CB}(X;M_{R\otimes R}(X\otimes
R\otimes R))$. Similarly, we can define {\em left coaction} as a
completely
bounded map $\gamma $ from $X$ to $M_{R}(R\otimes X)$ such that $(\mathrm{id}
\otimes \gamma )\gamma =(\delta \otimes \mathrm{id})\gamma $.

\smallskip
\noindent (b) A right coaction $\beta $ is said to be \emph{right}
(respectively, \emph{left}) \emph{non-degenerate} if the linear span
of $\{\beta (x)\cdot s:x\in X;s\in R\}$ (respectively, $\{s\cdot \beta
(x):x\in X;s\in R\}$) is norm dense in $X\otimes R$ (we recall that
$(l,r)\cdot s = l(s)$ and $s\cdot(l,r) = r(s)$ for any
$(l,r)\in M_R(X\otimes R)$ and $s\in R$).

\smallskip
\noindent (c) $X$ is said to be a
\emph{right $R$-comodule} (respectively, left $R$-comodule) if there exists a right coaction (respectively, left coaction) of $R$
on $X$.
\end{definition}

\medskip
The right coaction identity in part (a) actually means that $\Phi \circ
(\beta \otimes \mathrm{id})\circ \beta =(\mathrm{id}\otimes
\delta )\circ \beta $ (where $\Phi $ is the forgettable complete
isometry from $M_{R}(M_{R}(X\otimes R)\otimes R)$ to $M_{R\otimes
R}(X\otimes R\otimes R)$ given by Corollary \ref{1.11}(c)).

\medskip
\begin{lemma}
\label{A30} (a) The predual $\Re_{\ast }$ of $\Re$ has a
left (or right) identity if and only if it is unital.

\smallskip
\noindent (b) Suppose that $\beta $ is a right (or left) non-degenerate
right coaction of $S$ in $X$ and $\epsilon$ is a counit of $S$. Then $(
\mathrm{id}\otimes \epsilon)\circ \beta =\mathrm{id}$. The same is true for
a left coaction.

\smallskip\noindent
(c) Let $\beta $ be a right coaction of $S$ on $X$ and let $Y$ be a closed
subspace of $X$.
Suppose that the canonical quotient map $q$ from $X$ to
$X/Y$ satisfies the following condition:
$(q\otimes \id)\beta(Y) = (0)$.
Then $\beta$ induces a right coaction $\hat\beta$ of $S$ on $X/Y$ such that
$\hat\beta\circ q = (q\otimes\id)\circ \beta$.
\end{lemma}
\noindent \textbf{Proof:} (a) Suppose that $\epsilon$ is a left identity of $
\Re_{\ast }$.
Then we have $(\nu \otimes \epsilon)\delta ((\mathrm{id}\otimes \omega
)\delta (s))=(\nu \otimes (\epsilon\otimes \omega )\circ \delta )\delta (s)=(\nu
\otimes \omega )\delta (s)$ (for any $\omega ,\nu \in \Re_{\ast }$ and $
s\in \Re$).
Now weak*-density of $\delta (\Re)(1\otimes \Re)$ in
$\Re\bar{\otimes}\Re$ implies that
$\epsilon$ is an identity of $\Re_{\ast }$.

\smallskip
\noindent (b) Suppose that $\beta$ is right non-degenerate.
For any $x\in X$ and $s\in S$, the tensor $x\otimes s$ can be
approximated by sums of elements of the form $\beta (y)\cdot t$ ($y\in
X$; $t\in S$) and hence $x$ can be approximated by sums of elements of the
form $(\mathrm{id}\otimes g)\beta (y)$ ($g\in S^{\ast }$; note
that $M_{S}(X\otimes S)$ can be regarded as a subspace of some $M(\mathcal{K}
(H)\otimes S)$ by Proposition \ref{1.d}). The lemma now follows from the
fact that $(\mathrm{id}\otimes \epsilon)\beta ((\mathrm{id}\otimes
g)\beta (y)) = (\id\otimes(\epsilon\otimes g)\delta)\beta(y) =
(\mathrm{id}\otimes g)\beta (y)$. The other three
cases can be proved similarly.

\smallskip\noindent
(c) By Lemma \ref{A1}, there exists a map $\hat\beta$ that satisfies the
required equality.
It remains to check the right coaction identity.
In fact, $(\hat\beta\otimes \id)\hat\beta \circ q =
(\hat\beta\otimes \id)(q\otimes \id)\beta =
((q\otimes \id)\beta\otimes \id)\beta =
(q\otimes \id\otimes \id)(\id\otimes \delta)\beta =
(\id\otimes\delta)(q\otimes \id)\beta$
(note that as we are working with multipliers of an operator bimodule
instead of elements in a $C^*$-algebra, cautions needed
to be taken for each of the equalities above; in particular, the last
equality follows from Lemma \ref{1.12}(b)).

\medskip
\begin{example}
\label{B1} (a) Let $\Gamma $ be a discrete group.
The reduced group $C^*$-algebra $C^*_r(\Gamma)$ is a Hopf $C^*$-algebra
with coproduct given by $\delta(\lambda_r)=\lambda_r\otimes\lambda_r$ (where
$\lambda_r$ is the canonical image of $r\in \Gamma$ in $C^*_r(\Gamma)$).
If $\beta $ is a
non-degenerate coaction of $C_{r}^{\ast }(\Gamma )$ on a $C^{\ast }$-algebra
$A$ (in the sense of Definition \ref{Hopf}(c)), then $A$ can be decomposed as
$A=\overline{\bigoplus_{r\in \Gamma }A_{r}
}$ (see \cite[2.6]{Ng0}). Let $F$ be any subset of $\Gamma $ and $A_{F}=
\overline{\bigoplus_{r\in F}A_{r}}$. Then the restriction $\beta _{F}$
of $\beta $ on $A_{F}$ is a right coaction of $C_{r}^{\ast }(\Gamma )$ on $
A_{F}$. Moreover, it is not hard to see that this right coaction is also
(2-sided) non-degenerate.

\smallskip
\noindent (b) Suppose that $\beta $ is a right (or left)
non-degenerate right coaction of $C_{r}^{\ast }(\Gamma )$ on $X$. Let $
X_{r}=\{x\in X:\beta (x)=x\otimes \lambda _{r}\}$. Then $X_{r}=(
\mathrm{id}\otimes \varphi _{r})\beta (X)$ where $\varphi _{r}$ is the
functional on $C_{r}^{\ast }(\Gamma )$ satisfying $\varphi _{r}(\lambda
_{t})=\delta _{r,t}$ (here $\delta _{r,t}$ means the Kronecker delta) as
defined in \cite[\S 2]{Ng0}. Now by the right (respectively, left)
non-degeneracy of $\beta $, we have $X=\overline{\bigoplus_{r\in
\Gamma }X_{r}}$.

\smallskip
\noindent (c) Let $G$ be a locally compact group.
Then $C_0(G)$ is a Hopf $C^*$-algebra with a coproduct defined by
$\delta(f)(s,t) = f(st)$ (note that $M(C_0(G)\otimes C_0(G)) =
C_b(G\times G)$).
Right coactions of $C_{0}(G)$ on $X$ are in one to one
correspondence with completely bounded
actions of $G$ on $X$ in the following sense: an action $\alpha $ of
$G$ on $X$ is said to be \emph{completely bounded} if

\begin{enumerate}
\item[i.]  there is $\lambda >0$ such that $\sup \{\Vert (\mathrm{id}_{
\mathcal{K}}\otimes \alpha _{t})(\bar{x})\Vert :t\in G\}\leq \lambda \Vert
\bar{x}\Vert $ for any $\bar{x}\in \mathcal{K}\odot X$;

\item[ii.]  $\alpha _{\bullet }(x)$ is a continuous map from $G$ to $X$
for any fixed $x\in X$;
\end{enumerate}
\noindent (or equivalently, \noindent $\alpha $ induces a bounded continuous
action of $G$ on the subspace $\mathcal{K}\odot X$ of
$\mathcal{K}\otimes X$).
Moreover, if the right coaction is a
complete isometry, then condition (i) is replaced by the following condition:

\begin{enumerate}
\item[i'.]  $\sup \{\Vert (\mathrm{id}_{\mathcal{K}}\otimes \alpha _{t})(
\bar{x})\Vert :t\in G\}=\Vert \bar{x}\Vert $ for any $\bar{x}\in \mathcal{K}
\odot X$.
\end{enumerate}
\noindent Indeed, by considering $X$ as a closed subspace of a $C^{\ast }$-algebra, we see that $M_{C_{0}(G)}(X\otimes C_{0}(G))=C_{b}(G;X)$. Hence as
in the case of $C^{\ast }$-algebras, a right coaction $\delta $ induces an action $
\alpha $ of $G$ on $X$ such that $\alpha _{t}(x)=\delta (x)(t)$. Since $
\delta (X)\subseteq C_{b}(G;X)$ and $\delta $ is bounded, there exists $
\lambda ^{\prime }>0$ such that $\sup \{\Vert \alpha _{t}(x)\Vert :t\in
G\}\leq \lambda ^{\prime }\Vert x\Vert $ and for fixed $x\in X$, the map $\alpha
_{\bullet }(x)$ is continuous. As $\delta $ is
completely bounded, we can replace $X$ with $\mathcal{K}\otimes X$ and show that
$\alpha $ is a completely bounded action. Conversely, let $\alpha $ be a
completely bounded action. If we define $\delta (x)(t)=\alpha _{t}(x)$
(for any $x\in X$ and $t\in G$), then condition (ii) implies that $\delta
(x)\in C(G;X)$ (i.e. a continuous map) and so $k\otimes \delta (x)\in C(G;
\mathcal{K}\otimes X)$ (for any $k\in \mathcal{K}$). Furthermore, condition
(i) shows that $\mathrm{id}\otimes \delta $ is a bounded map from $\mathcal{K
}\odot X$ to $C_{b}(G;\mathcal{K}\otimes X)=M_{C_{0}(G)}(\mathcal{K}\otimes
X\otimes C_{0}(G))$. It is not hard to check that $\delta $ is a right coaction
and the correspondence is established. In this case, $\delta $ is injective
if and only if $\alpha _{t}$ is injective for all $t\in G$ (because
$\alpha _{s}\alpha _{t}=\alpha _{st}$) or 
equivalently, $\alpha_t$ is injective for some $t\in G$.
It is the case if and only if $\alpha _{e}=I_{X}$ (note that $\alpha
_{e}(\alpha _{e}(x)-x)=0$).

\smallskip
\noindent (d) Again when $S=C_{0}(G)$, there is an one to one correspondence
between right coactions of $S$ on $X$ and completely bounded representations of $G$
on $X$: a map $T$ from $G$ to $\mathrm{CB}(X;X)$ is called a
\emph{completely bounded representation} if

\noindent \hspace{1em} i. $T_{r}\circ T_{s}=T_{rs}$;

\noindent \hspace{1em} ii. $\sup \{\Vert T_{r}\Vert _{\mathrm{cb}}:r\in
G\}<\infty $;

\noindent \hspace{1em} iii. the map $T_{\bullet }(x)$
from $G$ to $X$ is continuous for any $x\in X$.

\noindent In particular, if $H$ is a Hilbert space and $H_{c}$ is the column
operator space of $H$, then right coactions of $C_{0}(G)$ on $H_{c}$ are exactly
bounded continuous representations of $G$ on $\mathcal{L}(H)=\mathrm{CB}
(H_{c};H_{c})$.

\smallskip\noindent
(e) By \cite[2.4]{Ng3}, a coaction of $S$
on a Hilbert $C^*$-modules $E$ (in the
sense of \cite{BS0}) defines a right coaction of $S$ on
the ``column space'' $E_c$.

\smallskip\noindent
(f) By \cite[2.5 \& 2.7]{Ng3}, if a corepresentation of $S$ on a Hilbert space $H$
is either unitary or ``non-degenerate'', then it gives rise to a coaction of $S$
on $H_c$.

\end{example}

\bigskip

\bigskip

\section{Cohomology of Hopf $C^{\ast }$-algebras}

\bigskip

In this section, we will define and study cohomology theories for Hopf $
C^{\ast }$-algebras with coefficients in their bicomodules.
Let $X$ be an operator space and let $\beta$ and $\gamma$ be respectively a
right and a left coactions of $S$ on $X$.
Then $(X,\beta, \gamma)$ is said to be a $S$\emph{-bicomodule} if
$(\mathrm{id}\otimes
\beta )\gamma =(\gamma \otimes \mathrm{id})\beta $.
Let us first consider a straight forward way to define cohomology
(which is a ``dual analogue'' of Banach algebra homology).

\bigskip
\noindent{\bf Notation:}
\emph{For simplicity, we may sometimes use $X$ to denote the bicomodule
$(X,\beta,\gamma)$.
Throughout this section, $S^{n}$ is the n-times spatial tensor
product of $S$ (whereas $S^{0}=
\mathbb{C}$) and $\sigma _{n,k}$ ($1\leq k\leq n$) is the completely bounded
map from $M_{S^{n}}(S^{k}\otimes X\otimes S^{n-k})$ to $M_{S^{n}}(X\otimes
S^{n})$ defined by $(s_{n-k+1}\otimes ...\otimes s_{n}\otimes x\otimes
s_{1}\otimes ...\otimes s_{n-k})^{\sigma _{n,k}}=x\otimes s_{1}\otimes
...\otimes s_{n}$ (see Lemma \ref{1.12}(a)).}

\bigskip

For $n\geq 1$, we define a
completely bounded map $\delta _{n}$ from $M_{S^{n}}(X\otimes S^{n})$ to $
M_{S^{n+1}}(X\otimes S^{n+1})$ by
\begin{eqnarray*}
\delta _{n}(x\otimes s_{1}\otimes...\otimes s_{n})
&=&
\beta (x)\otimes s_{1}\otimes ...\otimes s_{n}+\sum_{k=1}^{n}(-1)^{k}
x\otimes s_{1}\otimes ...\otimes s_{k-1}\otimes\delta (s_{k})\otimes s_{k+1}
\otimes ...\otimes s_{n}+\\
& &
\qquad \qquad
(-1)^{n+1}(\gamma(x)\otimes s_{1}\otimes ...
\otimes s_{n})^{\sigma _{n+1,1}}
\end{eqnarray*}
and $\delta_{0}(x)=\beta (x)-\gamma (x)^{\sigma _{1,1}}$
($\delta _{n}$ is well
defined by Lemma \ref{1.12}(a) and Corollary \ref{1.11}(c)). We need to show
that $(M_{S^{n}}(X\otimes S^{n}),\delta _{n})$ is a cochain complex.

\medskip
\begin{lemma}
\label{B8} $\delta _{n}\circ \delta _{n-1}=0$ for $n=1,2,3,...$.
\end{lemma}
\noindent \textbf{Proof:}
Note, first of all, that $\delta _{1}\circ \delta
_{0}(x)=(\beta \otimes \mathrm{id})\beta (x)-(\mathrm{id}\otimes
\delta )\beta (x)+(\gamma \otimes \mathrm{id})\beta (x)^{\sigma
_{2,1}}-(\mathrm{id}\otimes \beta )\gamma (x)^{\sigma _{2,1}}+(\delta
\otimes \mathrm{id})\gamma (x)^{\sigma _{2,2}}-(\gamma \otimes \mathrm{id}
)(\gamma (x)^{\sigma _{1,1}})^{\sigma _{2,1}}=0$.
This established the equality for $n=1$.
For the case of $n>1$, the crucial point is to show that
$\sum_{k=1}^{n}\sum_{i=1}^{n-1}(-1)^{k+i}(\mathrm{id}
_{X}\otimes \mathrm{id}^{k-1}\otimes \delta \otimes \mathrm{id}^{n-k})\circ (
\mathrm{id}_{X}\otimes \mathrm{id}^{i-1}\otimes \delta \otimes \mathrm{id}
^{n-i-1})=0$.
This can be doned by a decomposition (into a sum of
summations according to the relative positions of the $i$'s and $k$'s
in the original summation) as well as a tedious comparison.

\bigskip

Now we can define a cohomology
$H^{n}(S;X)=\mathrm{Ker}(\delta _{n})/\mathrm{Im}(\delta _{n-1})$ for
$n\in \mathbb{N}$ and
$H^{0}(S;X)=\{ x\in X:\beta (x)=\gamma (x)^{\sigma _{1,1}}\}$.
It is called the
\emph{natural cohomology }\emph{of} $S$ \emph{with coefficient in the
bicomodule }$(X,\beta ,\gamma )$.

\medskip
\begin{example}
\label{B17} (a) Let $G$ be a locally compact group and $S=C_{0}(G)$.
We have already seen in Example \ref{B1}(c) that $S$ is a Hopf $C^*$-algebra and 
a right $S$-comodule $X$ is a
``completely bounded left $G$-module'' and $M_{S}(X\otimes S)=C_{b}(G;X)$. Hence
a $S$-bicomodule is a ``completely bounded $G$-bimodule''. In this case, $\psi
\in C_{b}(G;X)$ is in $Ker(\delta _{1})$ if and only if it is a derivation
in the sense that $\psi (st)=s\cdot \psi (t)+\psi (s)\cdot t$ $(s,t\in G)$.
Moreover, $\psi \in Im(\delta _{0})$ if and only if there exists $x\in X$
such that $\psi (s)=s\cdot x-x\cdot s$ $(s\in G)$. Hence $H^{1}(S;X)$ is
formally a kind of group cohomology of $G$.

\smallskip
\noindent (b) If $S=c_{0}(\Gamma )$ for a discrete group $
\Gamma $ and the left coaction $\gamma$ is $1\otimes \mathrm{id}_{X}$,
then $H^{1}(c_{0}(\Gamma
);X)$ coincides with the first group cohomology $H^{1}(\Gamma ;X)$
(with coefficients in the $c_0(\Gamma)$-bicomodule $X$) that studied in
\cite{Evens} (see \cite[p.9]{Evens}).

\smallskip
\noindent (c) On the other hand, if $G$ is a profinite group and $\gamma
=1\otimes \mathrm{id}_{X}$, the two groups $H^{n}(C_{0}(G);X)$ and
$H^{n}(G;X)$ coincide where $H^{n}(G;X)$ is the group cohomology studied
in \cite[\S 9.1]{Wilson} and \cite[\S 2.2]{Serre}.
\end{example}

\medskip
It is natural to consider the cohomology of a ``dual bicomodule'' and expect
it to relate to the amenability of the Hopf $C^{\ast }$-algebra.
However, it seems impossible to define dual comodule structure on the dual
space of a $S$-bicomodule. Nevertheless, we still have a kind of ``dual
cohomology theory'' with the cochain complex starting with the dual space.

\medskip

Let $\partial _{n}$ be the completely bounded map from $\mathrm{CB}
(X;M(S^{n}))$ to $\mathrm{CB}(X;M(S^{n+1}))$ defined by
$$\partial_{n}(T)=
\cases{
(T\otimes \mathrm{id})\circ \beta + \sum_{k=1}^{n}(-1)^{k}
(\mathrm{id}^{n-k}\otimes \delta \otimes \mathrm{id}^{k-1})\circ T+(-1)^{n+1}
(\mathrm{id}\otimes T)\circ \gamma  & \quad $n=1,2,3...$\cr
(T\otimes id)\circ \beta -(\mathrm{id}\otimes T)\circ \gamma  & \quad $n=0$
}
$$
for any $T\in \mathrm{CB}(X;M(S^{n}))$. It is well defined by Lemma \ref
{1.12}(a) and Corollary \ref{1.11}(c). Again, it gives a cochain complex.

\medskip
\begin{lemma}
\label{B10} $\partial _{n}\circ \partial _{n-1}=0$ for $n=1,2,3,...$.
\end{lemma}

\medskip
We can now define another cohomology theory by $H_{d}^{n}(S;X)=\mathrm{Ker}
(\partial _{n})/\mathrm{Im}(\partial _{n-1})$ ($n=1,2,3,...)$ and $
H_{d}^{0}(S;X)=\{f\in X^{\ast }:(f\otimes \mathrm{id})\circ \beta =(
\mathrm{id\otimes }f)\circ \gamma \}$. It is called the \emph{dual
cohomology }\emph{of }$S$ \emph{with coefficient in }$(X,\beta ,\gamma
)$. The name comes from the fact that it is a cohomology theory ``with
coefficient in
the dual space $X^{\ast }$''. Furthermore, the Hopf von Neumann algebra
analogue of this cohomology can actually be regarded as a dual cohomology
theory (see
Proposition \ref{C10} and Proposition \ref{C15} below).

\medskip

Note that the idea of the dual cohomology is similar to the cohomology theory
of coalgebras studied in \cite{Doi} but although they look alike, they
``behave differently'' even in the case of discrete groups
(see Remark \ref{rmstrcd}(c)).

\medskip
\begin{example}
\label{B18} Let $\Gamma $ be a discrete group.
Suppose that $(X,\beta ,\gamma )$ is a $C_{r}^{\ast }(\Gamma )$-bicomodule
(see Example \ref{B1}(a)\&(b)) such that
both $\beta $ and $\gamma $ are either left or right non-degenerate.
Then by Example \ref{B1}(b), $X$ can be decomposed into two directed sums $
\overline{\oplus _{s\in \Gamma }X_{s}^{\beta }}$ and $\overline{\oplus
_{t\in \Gamma }X_{t}^{\gamma }}$ corresponding to $\beta $ and $\gamma
$ respectively (where $X_{s}^{\beta }$=(\textrm{id}$\otimes \varphi
_{s})\beta (X)$ and $\varphi _{s}$ is the functional as defined in Example
\ref{B1}(b)). Moreover, since $(\mathrm{id}\otimes \beta )\gamma
=(\gamma \otimes \mathrm{id})\beta $, for any $s,t\in \Gamma $, the spaces
$X_{s}^{\beta }$ and $X_{t}^{\gamma }$ can be decomposed further into $
\overline{\oplus _{r\in \Gamma }(X_{s}^{\beta })_{r}}$ and $\overline{
\oplus _{r\in \Gamma }(X_{t}^{\gamma })_{r}}$ corresponding to $\gamma $ and
$\beta $ respectively such that $(X_{s}^{\beta
})_{t}=(X_{t}^{\gamma })_{s}$. For any $x\in (X_{s}^{\beta })_{t}$, we have
$\beta (x)=x\otimes \lambda _{s}$ and $\gamma (x)=\lambda _{t}\otimes x$. 
Let $\alpha \in Ker(\partial _{1})$. Then
\[
\delta (\alpha (x))=\alpha (x)\otimes \lambda _{s}+\lambda _{t}\otimes
\alpha (x)
\]
and so $(\varphi _{r}\otimes \mathrm{id})\delta (\alpha (x))=
\varphi _{r}(\alpha (x))\lambda _{s}+\delta_{r,t}\alpha (x)$
for any $r\in \Gamma$ (where $\delta _{r,t}$ is the Kronecker delta).
Thus by putting $r=t$, we obtain $\alpha (x)=\varphi _{t}(\alpha (x))(\lambda _{t}-
\lambda_{s})$ (note that $(\varphi_{t}\otimes\mathrm{id})\delta (\alpha(x))=
\varphi _{t}(\alpha(x))\lambda _{t}$).
Similarly, we have $\alpha (x)=\varphi _{s}(\alpha (x))(\lambda
_{s}-\lambda _{t})$. Now for any $y\in \oplus _{s\in \Gamma
}X_{s}^{\beta }$, we define
\[
f(y)=\sum \varphi _{s}(\alpha (y_{s}))
\]
(where $y=\sum y_{s}$ and $y_{s}\in X_{s}^{\beta }$). If $f$ extends
to a continuous function on $X$, then it is not hard to see that $\partial
_{0}(f)=\alpha $. In fact, for any $x\in (X_{t}^{\gamma })_{s}$, we have
$(\id\otimes f)\gamma (x)=f(x)\lambda _{t}=\varphi _{s}(\alpha (x))\lambda _{t}$
(by the definition of $f$ and the fact that $x\in X^\beta_s$) and ($f\otimes
\mathrm{id})\beta (x)=f(x)\lambda _{s}=\varphi _{s}(\alpha (x))\lambda
_{s}$ (as $(X_{t}^{\gamma })_{s}=(X_{s}^{\beta })_{t}$). It is not
clear for the moment whether all such functions defined in this way are
continuous (this means that $H^1_d(C^*_r(\Gamma;X)=(0)$).
However, we will see in Corollary \ref{discrete} that this will imply the
amenability of $\Gamma $.
The converse is also true because if $f\in X^*$ such that
$\partial_0(f)=\alpha$, then
$f(x)=\varphi_s(\alpha(x))$ for any $x\in (X^\beta_s)_t$.
\end{example}

\medskip
Next, we would like to study the situation when these two cohomology
theories vanish. First of all, we will consider the case when the left
coaction $\gamma $ of the $S$-bicomodule $(X,\beta ,\gamma )$ is
trivial in the sense that $\gamma =0$ (there is another meaning for
the triviality: $\gamma =1\otimes \mathrm{id}_{X}$ but we will
not consider this situation until Section 4).
In this case, the corresponding (one-sided) natural and dual
cohomologies will be
denoted by $H_{r}^{n}(S;X)$ and $H_{d,r}^{n}(S;X)$ respectively.
We have the following simple result concerning these one-sided
cohomologies.

\medskip
\begin{proposition}
\label{B12} Let $(S,\delta)$ be a saturated Hopf $C^*$-algebra and
$(X,\beta )$ be a right $S$-comodule.

\smallskip
\noindent (a) $H_{r}^{0}(S;X)=(0)$ if and only if $\beta $ is
injective. Moreover, if $\beta $ is either left or right
non-degenerate, then $H_{d,r}^{0}(S;X)=(0)$.

\smallskip
\noindent (b) If $(S,\delta)$ is counital, then $H_{r}^{n}(S;X)=(0)$ and $
H_{d,r}^{n}(S;X)=(0)$ ($n\geq 1)$.

\smallskip
\noindent (c) If $H_{d,r}^{1}(S;S)=(0)$, then $(S,\delta)$ is counital.
\end{proposition}
\noindent
\textbf{Proof:} (a) The first statement is obvious. To show the
second statement, we suppose that $\beta$ is right non-degenerate
and take any $f\in \mathrm{Ker}(\partial _{0})$ (i.e. $
(f\otimes \mathrm{id})\beta =0)$.
Now the density of $\beta (X)\cdot S$ in
$X\otimes S$ will imply that $f$ is zero.
The argument for the case when $\beta$ is left non-degenerate is similar.

\smallskip
\noindent (b) Suppose that $\epsilon$ is the counit of $S$. Let $m\in
M_{S^{n}}(X\otimes S^{n})$ be such that $\delta _{n}(m)=0$ and let $m_{\epsilon}=(
\mathrm{id}_{X}\otimes \mathrm{id}^{n-1}\otimes \epsilon)(m)$.
Then it is clear by Lemma \ref{1.12}(b) that
\[
(\beta \otimes \mathrm{id}^{n-1})(m_{\epsilon})=(\mathrm{id}_{X}\otimes
\mathrm{id}^{n}\otimes \epsilon)(\beta \otimes \mathrm{id}
^{n})(m)=(-1)^{n-1}m-\sum_{k=1}^{n-1}(-1)^{k}(\mathrm{id}_{X}\otimes \mathrm{
id}^{k-1}\otimes \delta \otimes \mathrm{id}^{n-k-1})(m_{\epsilon}).
\]
Hence $m=\delta _{n-1}((-1)^{n-1}m_{\epsilon})$. This shows that $
H_{r}^{n}(S;X)=(0)$. On the other hand, let $T\in \mathrm{CB}(X;M(S^{n}))$
such that $\partial _{n}(T)=0$. If $F=(\epsilon\otimes \mathrm{id}^{n-1})\circ T$,
then it is clear from $(\epsilon\otimes \mathrm{id}^{n})(\partial _{n}(T))=0$ that
\[
(F\otimes \mathrm{id})\circ \beta +(-1)^{n}T+\sum_{k=1}^{n-1}(-1)^{k}(
\mathrm{id}^{n-k-1}\otimes \delta \otimes \mathrm{id}^{k-1})\circ F=0.
\]
Hence $T=(-1)^{n-1}\partial _{n-1}(F)$.

\smallskip
\noindent (c) Consider $\mathrm{id}\in \mathrm{CB}(S;M(S))$. Then clearly $
\partial _{1}(\mathrm{id})=0$ and so there exists $\epsilon\in S^{\ast }$ such that
$(\epsilon\otimes \mathrm{id})\circ \delta =\mathrm{id}$. This shows that $\epsilon$ is a
left identity for $S^{\ast }$ and hence a two-sided identity (by Lemma \ref
{A30}(a) and the fact that $S^{**}$ is a saturated Hopf von Neumann algebra).

\bigskip

In the following, we will study 2-sided dual cohomology. For the moment, we
only have the complete picture for the case when $S$ is unital (i.e. it
represents a compact quantum group) and a partial
picture if $S$ has property (S) (in the sense
of \cite{Was}).
In these cases, the vanishing of the dual cohomology is
related to the existence of codiagonals defined as follows.

\medskip
\begin{definition}
\label{B19} Suppose that $R$ is a Hopf $C^*$-algebra with counit $\epsilon$.
Let $Y$ be a subspace of $M(R\otimes R)$ containing $\delta (R)$ such
that $(\id\otimes(\id\otimes f)\circ \delta)(Y)\subseteq Y$ and
$((f\otimes \id)\circ \delta\otimes \id)(Y)\subseteq Y$ for
any $f\in R^*$.
Then $F\in Y^{\ast }$
is said to be a \emph{codiagonal} if $F\circ \delta =\epsilon$
on $R$ and $F\circ (\mathrm{id}\otimes (\mathrm{id}\otimes g)\circ
\delta )=F\circ ((g\otimes \mathrm{id})\circ \delta
\otimes \mathrm{id})$ on $Y$ (for any $g\in R^{\ast })$.
\end{definition}

\medskip
We will defer the illustrations and examples for the
codiagonals until Example \ref{B23} and Remark \ref{rmstrcd}(a) \& (b).

\medskip

Note that from this point on, we will need quite a lot of materials
from the appendix (unless the readers want to confine themselves in
the case of unital Hopf $C^*$-algebras -- in which case, please see
part (b) and (c) of the following Remark).
Therefore, perhaps it will be a good idea if the
readers can digress to the Appendix at this point (we are sorry that since the
materials in the appendix are a bit technical and not in the same
favour as the other parts of this paper, we decided to 
study them in the appendix).

\medskip
\begin{remark}
\label{ext-Y}
Let $Y=\hat{U}(R\otimes R)$ (see Remark \ref{U(A)}).

\smallskip\noindent
(a) If $R$ has property (S) (in particular, if $R$ is a nuclear $C^{\ast }$-algebra),
then by Lemma \ref{2sU(A)},
$Y$ is the biggest unital $C^*$-subalgebra
of $M(R\otimes R)$ for which $\id\otimes \delta$ and $\delta\otimes\id$
can be extended.

\smallskip
\noindent
(b) If $R$ is unital, then $Y=R\otimes R$ and $\delta\otimes\id$ and
$\id\otimes\delta$ obviously define on $Y$ (without $R$ having
the property (S)).

\smallskip
\noindent
(c) In both of the cases (a) and (b) above, $Y$ satisfies the conidtion in
Definition \ref{B19} and $F\in Y^{\ast }$ is
a codiagonal if and only if $F\circ \delta =\epsilon$ and $
(F\otimes \mathrm{id})\circ (\mathrm{id}\otimes \delta )=(\mathrm{id}\otimes
F)\circ (\delta \otimes \mathrm{id})$.
\end{remark}

\medskip
\begin{proposition}
\label{exist-sc} Let $S$ be a saturated Hopf $C^*$-algebra. Suppose that $S$
either is unital or has property (S).

\smallskip
\noindent (a) If $H_{d}^{1}(S;X)=(0)$ for any $S$-bicomodule $X$, then there
exist a counit $\epsilon$ on $S$ as well as a codiagonal $F$ on $\hat{U}
(S\otimes S)$.
Moreover, if $S$ is unital, we obtain the same conclusion
even if $H_d^1(S;X)$ vanish only for those $S$-bicomodules
$(X,\beta,\gamma)$ such that both $\beta$
and $\gamma$ are 2-sided non-degenerate.

\smallskip
\noindent (b) If there exist a counit $\epsilon$ on $S$ and a codiagonal 
$\bar{F}$ on $M(S\otimes S)$ such that $\bar F\circ \delta = \epsilon$ on $M(S)$, 
then $H_{d}^{1}(S;X)=(0)$ for any $S$-bicomodule
$(X,\beta ,\gamma )$ such that either $\beta $ or $
\gamma $ is left (or right) non-degenerate.
\end{proposition}
\noindent \textbf{Proof:} (a) By Proposition \ref{B12}(c), $S$ has a counit $
\epsilon$. Let $U(S)$ be the space $U_{\delta ,\delta }(S)$
(see Remark \ref{U(A)}).
It is not hard to see that $\delta (U(S))\subseteq \hat{U}
(S\otimes S)$ (since $(\mathrm{id}\otimes \delta )\delta (m)(1\otimes
1\otimes s)=(\delta \otimes \mathrm{id})(\delta (m)(1\otimes s))\in
\delta(U(S))\otimes S$ and
$(\delta \otimes \mathrm{id})\delta (m)(s\otimes 1\otimes 1)=(\mathrm{id}
\otimes \delta )(\delta (m)(s\otimes 1))\in S\otimes \delta(U(S))$
for any $m\in U(S)$).
Let $X$ be the quotient $\hat{U}(S\otimes S)/\delta (U(S))$ with the canonical quotient map $q$.
By Remark \ref{ext-Y}, $\mathrm{id}\otimes \delta $
induces a right coaction on $\hat{U}(S\otimes S)$.
Using the first equality
above, we have $(q\otimes \mathrm{id})(\mathrm{id}\otimes \delta )\delta
(U(S))=(0)$ in $M_{S}(X\otimes S)$ and by Lemma \ref{A30}(c), $\mathrm{id}
\otimes \delta $ induces a right coaction $\beta $ on $X$. Similarly, $
\delta \otimes \mathrm{id}$ induces a left coaction $\gamma$
on $X$. It is clear that
$(X,\beta ,\gamma )$ is a $S$-bicomodule. Consider the completely
bounded map $T=\epsilon\otimes \mathrm{id}-\mathrm{id}\otimes \epsilon$ from $\hat{U}
(S\otimes S)$ to $M(S)$. Since $T\circ \delta =0$, it induces a map $\hat{T}
\in \mathrm{CB}(X;M(S))$ (by Lemma \ref{A1}). Now for any $m\in \hat{U}
(S\otimes S)$,
\[
\partial _{1}(\hat{T})(q(m))=(T\otimes \mathrm{id})(\mathrm{id}\otimes
\delta )(m)-\delta ((\epsilon\otimes \mathrm{id})(m)-(\mathrm{id}\otimes \epsilon)(m))+(
\mathrm{id}\otimes T)(\delta \otimes \mathrm{id})(m)=0.
\]
Hence there exists $G\in X^{\ast }$ such that $\partial _{0}(G)=\hat{T}$
i.e. $\epsilon\otimes \mathrm{id}-\mathrm{id}\otimes \epsilon=(G\circ q\otimes \mathrm{id}
)(\mathrm{id}\otimes \delta )-(\mathrm{id}\otimes G\circ q)(\delta \otimes
\mathrm{id})$. Let $F=\epsilon\otimes \epsilon-G\circ q\in \hat{U}(S\otimes S)^{\ast }$.
Then
\[
(F\otimes \mathrm{id})(\mathrm{id}\otimes \delta )=\epsilon\otimes \mathrm{id}
-(G\circ q\otimes \mathrm{id})(\mathrm{id}\otimes \delta )=\mathrm{id}
\otimes \epsilon-(\mathrm{id}\otimes G\circ q)(\delta \otimes \mathrm{id})=(\mathrm{
id}\otimes F)(\delta \otimes \mathrm{id})
\]
\medskip and $F\circ \delta =(\epsilon\otimes \epsilon)\circ \delta -G\circ q\circ \delta
=\epsilon$.
It is easy to see that if $S$ is unital, then $\beta$ and $\gamma$
in the above are left as well as right non-degenerate.

\smallskip
\noindent (b) Suppose that $T\in \mathrm{CB}(X;M(S))$ such that $\partial
_{1}(T)=0$, i.e. $\delta \circ T=(T\otimes \mathrm{id}_{S})\circ \beta
+(\mathrm{id}_{S}\otimes T)\circ \gamma $. Let $\beta $ be left
non-degenerate and $f=\bar{F}\circ (T\otimes \mathrm{id})\circ \beta
\in X^{\ast }$. Then by the properties of $\beta $ and $\gamma $ as
well as the definition of codiagonal, for any $g\in S^{\ast }$ and $
x\in X$,
\begin{eqnarray*}
g(\partial _{0}(f)(x)) &=& \bar{F}((\mathrm{id}_{S}^{2}\otimes g)
((T\otimes \mathrm{id}_{S})\beta\otimes \id_S)\beta (x))-\bar{F}
((g\otimes \mathrm{id}_{S}^{2})(\mathrm{id}_{S}\otimes T\otimes \mathrm{id}
_{S})(\mathrm{id}_{S}\otimes \beta )\gamma (x)) \\
&=&\bar{F}((\mathrm{id}_{S}^{2}\otimes g)(\mathrm{id}
_{S}\otimes \delta )(T\otimes \mathrm{id}_{S})\beta (x))-\bar{F}
((g\otimes \mathrm{id}_{S}^{2})(\mathrm{id}_{S}\otimes T\otimes \mathrm{id}
_{S})(\mathrm{id}_{S}\otimes \beta )\gamma (x)) \\
&=&\bar{F}((g\otimes \mathrm{id}_{S}^{2})(\delta \otimes \mathrm{id}
_{S})(T\otimes \mathrm{id}_{S})\beta (x))-\bar{F}((g\otimes \mathrm{id}
_{S}^{2})(\mathrm{id}_{S}\otimes T\otimes \mathrm{id}_{S})(\gamma \otimes
\mathrm{id}_{S})\beta (x)) \\
&=&\bar{F}\circ (g\otimes \mathrm{id}_{S}^{2})[(T\otimes \mathrm{id}
_{S}^{2})(\beta \otimes \mathrm{id}_{S})\beta (x)+(\mathrm{id}
_{S}\otimes T\otimes \mathrm{id}_{S})(\gamma \otimes \mathrm{id}
_{S})\beta (x)]- \\
&&\qquad \bar{F}((g\otimes \mathrm{id}_{S}^{2})(\mathrm{id}_{S}\otimes
T\otimes \mathrm{id}_{S})(\gamma \otimes \mathrm{id}_{S})\beta (x)) \\
&=&\bar{F}(\delta ((g\otimes \mathrm{id}_{S})(T\otimes \mathrm{id}
_{S})\beta (x)))\quad =\quad g(T(\mathrm{id}_{X}\otimes \epsilon)\beta
(x))\quad =\quad g(T(x)).
\end{eqnarray*}
\noindent The final equality follows from Lemma \ref{A30}(b). The case when $\gamma $
is left non-degenerate can be proved similarly by using $f=\bar{F}\circ (
\mathrm{id}\otimes T)\circ \gamma $.

\bigskip

Note that the proof of part (a) is similar to that of the dual situation
for the existence of diagonals (see e.g. \cite{Pier}).

\medskip
\begin{remark} This proposition applies in particular to the case when
$S$ is the Hopf $C^*$-algebra $S_V$ defined in \cite[1.5]{BS}
(see also \cite[3.8]{BS}) for a coamenable regular multiplicative
unitary $V$ (in this case, $S=S_V$ is a nuclear $C^{\ast }$-algebra
by \cite{Ng}).
Note that this includes the situation of $C_0(G)$ for a locally compact
group $G$ (which need not be amenable).
\end{remark}

\medskip

In the case when the Hopf $C^*$-algebra $S$ is unital, we call 
a net $\{ F_i\}$ in $(S\otimes S)^*$ a 
{\it bounded approximate codiagonal} if $\{ \| F_i \|\}$ is bounded and for any 
$f\in S^*$, both $\|(f\otimes F_i)\otimes(\delta\otimes \id) - (F_i \otimes 
f)\circ (\id\otimes \delta)\|$ and $\|(F_i\otimes f)\circ (\delta\otimes 
\id)\circ \delta - f\|$ converge to zero.

\medskip
\begin{theorem}
\label{B20} Suppose that $(S,\delta)$ is a saturated unital Hopf
$C^{\ast }$-algebra. The following conditions are equivalent.

\begin{enumerate}
\item[(i)]  $S$ has a counit and $S\otimes S$ has a codiagonal;

\item[(ii)]  $S\otimes S$ has a bounded approximate codiagonal;

\item[(iii)]  For any $S$-bicomodule $(X,\beta ,\gamma )$ such that
either $\beta $ or $\gamma $ is right (or left) non-degenerate, $
H_{d}^{n}(S;X)=(0)$ for $n\geq 1$.
\end{enumerate}
\end{theorem}
\noindent \textbf{Proof:} It is clear that (i) implies (ii). To show that
(ii) implies (i), let $\left\{ F_{i}\right\} _{i\in I}$ be a bounded
approximate codiagonal of $(S\otimes S)^{\ast }$ and $F$ be a $
\sigma ((S\otimes S)^{\ast },S\otimes S)$-limit point of $\left\{
F_{i}\right\} _{i\in I}$.
By considering a subnet if necessary, we may assume that $
\left\{ F_{i}\right\} _{i\in I}$ $\sigma ((S\otimes S)^{\ast },S\otimes S)$-converges to $F$. For any $f\in S^{\ast }$ and $s,t\in S$, it
is easy to check that $(F\otimes f)(\mathrm{id}\otimes \delta )(s\otimes
t)=(f\otimes F)(\delta \otimes \mathrm{id})(s\otimes t)$ (note that $\delta
(s),\delta (t)\in S\otimes S$).\ Moreover, as $\delta ^{\ast }$ is $\sigma
((S\otimes S)^{\ast },S\otimes S)$-$\sigma (S^{\ast },S)$-continuous, for any $s\in S$ and $f\in S^{\ast }$, we have $(\delta ^{\ast
}(F)\otimes f)\delta (s)=\lim_{i}\delta ^{\ast
}(F_{i})((\mathrm{id}\otimes f)\delta (s))=\lim_{i}(m(F_{i})\cdot
f)(s)=f(s)$. Thus, $\delta ^{\ast }(F)$ is a left identity of $S^{\ast
}$ and hence a two-sided identity (by Lemma \ref{A30}(a)). 
By Proposition \ref{exist-sc}(a), we have that (iii) implies (i).
It remains to show that (i) implies (iii). In fact, the argument is
similar to that of Proposition \ref{exist-sc}(b).
Let $\epsilon$ be a counit on $S$ and $F$ be a codiagonal on
$S\otimes S$.
Suppose that $\beta $ is either left or right non-degenerate. For
any $T\in Ker(\partial
_{n})\subseteq
\mathrm{CB}(X;M(S^{n}))$, let $R=(\mathrm{id}^{n-1}\otimes F)\circ (T\otimes
\mathrm{id})\circ \beta $. Then,
\begin{eqnarray*}
\partial _{n-1}R &=&(\mathrm{id}^{n-1}\otimes F\otimes \id)(T\otimes \id^{2})(
\mathrm{id}\otimes \delta )\beta +\sum_{j=1}^{n-1}(-1)^{j}(\mathrm{id}
^{n-j-1}\otimes \delta \otimes \mathrm{id}^{j-1})(\mathrm{id}^{n-1}\otimes
F)(T\otimes \mathrm{id})\beta + \\
&&\qquad (-1)^{n}(\mathrm{id}^{n}\otimes F)(\mathrm{id}\otimes T\otimes
\mathrm{id})(\mathrm{id}\otimes \beta )\gamma  \\
&=&(\mathrm{id}^{n}\otimes F)(\mathrm{id}^{n-1}\otimes \delta \otimes
\mathrm{id})(T\otimes \mathrm{id})\beta +(\mathrm{id}^{n}\otimes
F)[\sum_{j=1}^{n-1}(-1)^{j}(\mathrm{id}^{n-j-1}\otimes \delta \otimes
\mathrm{id}^{j+1})(T\otimes \mathrm{id})\beta ]+ \\
&&\qquad (-1)^{n}(\mathrm{id}^{n}\otimes F)(\mathrm{id}\otimes T\otimes
\mathrm{id})(\gamma \otimes \mathrm{id})\beta  \\
&=&(\mathrm{id}^{n}\otimes F)[((T\otimes \mathrm{id})\beta +
(-1)^{n+1}(\mathrm{id}\otimes T)\gamma) \otimes \mathrm{id}]\beta +
(-1)^{n}(\mathrm{id}^{n}\otimes F)(\mathrm{id}\otimes T\otimes
\mathrm{id})(\gamma \otimes \mathrm{id})\beta  \\
&=&(\mathrm{id}^{n}\otimes F)(T\otimes \delta )\beta \quad =\quad
(T\otimes \epsilon)\beta \quad =\quad T.
\end{eqnarray*}
\noindent In the case when $\gamma $ is left (or right) non-degenerate, we
should instead use $R=(F\otimes \mathrm{id}^{n-1})\circ (\mathrm{id\otimes
T)\circ \gamma }$ in the above argument.

\bigskip

Using this result, we can show that all the dual cohomologies of a
Woronowicz AF algebra (see \cite[\S3]{Wang}) vanish 
(this will be proved in \cite[3.9]{Ng4}). 

\medskip
\begin{lemexp}
\label{B23} Let $\Gamma $ be any discrete amenable group.
Then $C^{\ast }(\Gamma)\otimes C^{\ast }(\Gamma )$
has a codiagonal.
\end{lemexp}
\noindent \textbf{Proof:}
Note that as $\Gamma$ is amenable,
$(C^{\ast }(\Gamma)\otimes C^{\ast }(\Gamma ))^*=B(\Gamma\times
\Gamma)$ (the
Fourier-Stieltjes algebra of $\Gamma \times \Gamma )$.
Let $F\in B(\Gamma \times \Gamma )$ and $f\in B(\Gamma )$.
Then $F\cdot f=f\cdot F$ if and only if $F(r,s)f(s)=f(r)F(r,s)$
for any $r,s\in \Gamma$.
By taking $f=\varphi_{r}$ (where $\varphi_r$ is as defined in 
Example \ref{B1}(b)), we see that the above equality is equivalent 
to $F(r,s)=0$ if $r\neq s$.
Moreover, $F\circ \delta =\epsilon$ if and only if $F(r,r)=1$. 
Now consider $F_{0}\in l_{\infty }(\Gamma \times \Gamma )$ defined by 
$F_{0}(r,s)=\delta_{r,s}$ (where $\delta_{r,s}$ is the Kronecker delta). 
It is not hard to see that $F_{0}$ is positive definite. Indeed,
suppose that $\left\{ (r_{1},s_{1}),...,(r_{n},s_{n})\right\} $ is any
finite set in $\Gamma \times \Gamma $. Then $
F_{0}((r_{i},s_{i})^{-1}(r_{j},s_{j}))=1$ if and only if $
r_{i}^{-1}r_{j}=s_{i}^{-1}s_{j}$. Define an equivalent relation $\sim $ on $
\left\{ (r_{1},s_{1}),...,(r_{n},s_{n})\right\} $ by $(r,s)\sim (u,v)$
whenever s$r^{-1}=vu^{-1}$. Then $\{F_{0}((r_{i},s_{i})^{-1}(r_{j},s_{j}))
\}_{i,j=1,...,n}$ is equivalent to a direct sum of square matrices having $1$
in all their entries.
Hence $F_{0}\in B(\Gamma \times \Gamma )_{+}$ and
$C^{\ast }(\Gamma )\otimes C^{\ast }(\Gamma )$ have a
codiagonal.

\medskip
\begin{remark}
\label{rmstrcd}(a) By Corollary \ref{cor-sc} in Section 4, for any amenable
group $G$, the space $\hat{U}(G\times G)$ has a codiagonal. Hence, $C(G)\otimes
C(G)$ has a codiagonal for any compact group $G$.

\smallskip
\noindent (b) Suppose that $G$ is a locally compact
amenable group.
The same argument as in the above lemma shows that if
$F\in B(G\times G)$ satisfying the second condition of
Definition \ref{B19}, then $F(r,s)=0$ if $r\neq s$.
Therefore, it seems inappropriate to consider codiagonal
on $S\otimes S$ (instead of $\hat U(S\otimes S)$)
for a general Hopf $C^{\ast }$-algebras $S$.

\smallskip
\noindent (c) After we finished this manuscript, we
discovered that an analogue of the dual cohomology for coalgebras
has already been studied in \cite[\S 3.1]{Doi} and
a similar equivalence between
condition (i) and condition (iii) of Theorem \ref{B20}
(i.e. \cite[Thm 3]{Doi}) was obtained in the
purely algebraic setting (but with a different proof).
However, by the argument of Lemma \ref{B23} and \cite[Thm 3]{Doi}, for any
discrete group $\Gamma $, all the cohomologies of the coalgebra
$l^{1}(\Gamma )$ considered in \cite{Doi} vanish
(note that the functional $F_0$ in the above Lemma is well defined in
$B(\Gamma\times\Gamma) = (C^{\ast }(\Gamma)\otimes_{\max}
C^{\ast }(\Gamma ))^*$ even if $\Gamma$ is not amenable and
the restriction of $F_0$ on $l^1(\Gamma)\odot l^1(\Gamma)$ is the functional
required in $\cite[Thm 3]{Doi}$).
(A direct proof for this vanishing statement can also be obtained by using
a similar argument as in Example \ref{B18} in which case we don't
need the function $f$ to be continuous).
This, together with the following corollary, shows that the
``dual cohomology theory'' for Hopf algebras behaves very
differently from the one for Hopf $C^{\ast }$-algebras.
\end{remark}

\medskip
\begin{corollary}
\label{discrete} Let $\Gamma $ be a discrete group. Then $\Gamma $ is
amenable if and only if $H_{d}^{n}(C_{r}^{\ast }(\Gamma );X)=(0)$ for any $
n\geq 1$ and any left (or right) non-degenerate
$C_{r}^{\ast }(\Gamma )$-bicomodule $X$ and
equivalently, $H_{d}^{1}(C_{r}^{\ast }(\Gamma );X)=(0)$ for any left (or right)
non-degenerate $C_{r}^{\ast }(\Gamma )$-bicomodule $X$.
\end{corollary}
\noindent \textbf{Proof: }If $\Gamma $ is amenable, then Lemma \ref{B23} and
Theorem \ref{B20} show that all the dual cohomology of $C_{r}^{\ast }
(\Gamma)$ vanish.
Now if $H^1_d(C^*_r(\Gamma);X)$ vanishes for
any left (or right) non-degenerate $C_{r}^{\ast }(\Gamma )$-bicomodule
$(X,\beta,\gamma)$ (and in particular, when $\gamma=0$), then
$\Gamma$ is amenable by Proposition \ref{B12}(c).

\bigskip

Part of the above corollary (more precisely, the case when $n=1$)
is true for general locally compact groups (see
Theorem \ref{D9}(b)). Moreover, by Theorem \ref{D9}(a), the amenability of $
G $ is also equivalent to the vanishing of $H_{d}^{1}(C_{0}(G);X)$.

\medskip
\begin{remark}
Note that if $\Gamma$ is a discrete group such that $C^*(\Gamma)\otimes
C^*(\Gamma) = C^*(\Gamma)\otimes_{\max} C^*(\Gamma)$, then the argument
in Lemma \ref{B23} also gives the existence of a codiagonal on
$C^*(\Gamma)\otimes C^*(\Gamma)$ and in this case, the dual cohomologies
of $C^*(\Gamma)$ vanish.
Hence, the vanishing of the dual cohomologies of
$C^*(\Gamma)$ seems not strong enough to ensure the amenability of
$\Gamma$.
\end{remark}

We end this section with the following natural question: \emph{in general,
is there any relation between the vanishing of the dual cohomology and the
amenability or coamenability (see \cite{Ng2}) of the Hopf }$C^{\ast }$\emph{-algebra?} Some partial answers will be given in Section 4.

\bigskip

\bigskip

\section{Coactions and cohomology of Hopf von Neumann algebras}

\bigskip

In this section, we will study coactions and cohomology theories of
Hopf von Neumann algebras. We begin with coactions on dual operator
spaces (which is a natural generalisation of ordinary coactions on von
Neumann algebras).

\bigskip
\noindent
{\bf Notation:}
{\em Throughout this section, $\mathcal{X}$ is the dual
operator space of an operator space $\mathcal{X}_{\ast }$
and we recall from section 1 that $(\Re, \delta)$
is a Hopf von Neumann algebra.}

\bigskip
In order to define a coaction, we need to decide first of all, the range of
it (as in the case of Hopf $C^{\ast }$-algebras). Note that the range of a
coaction on a von Neumann algebra $\mathcal{M}$ by $\Re$ is the von
Neumann algebra tensor product $\mathcal{M}\overline{\otimes }\Re$.
For dual operator spaces, we have the following generalisation.
By \cite[2.1]{Blech}, there exists a weak*-homeomorphic complete
isometry from $\mathcal{X}$ to some $\mathcal{L}(H)$.
Let $\Re$ be represented as a von Neumann subalgebra of $\mathcal{L}(K)$
and let $\mathcal{X}\otimes _{\mathcal{F}}\Re=\{\alpha \in \mathcal{L}
(H)\overline{\otimes }\mathcal{L}(K):(\mathrm{id}\otimes \omega )(\alpha
)\in \mathcal{X};(\nu \otimes \mathrm{id})(\alpha )\in \Re$ for any $
\omega \in \Re_{\ast }$ and $\nu \in \mathcal{X}_{\ast }\}$ be the
\emph{Fubini product} (see \cite[p.188]{Ruan}). We recall from \cite[3.3]
{Ruan} that the Fubini product is independent of the representations (of $
\mathcal{X}$ and $\Re)$ and is the dual space of $\mathcal{X}_{\ast }
\hat{\otimes}\Re_{\ast }$. Moreover, if $\mathcal{X}$ is a von
Neumann algebra, then $\mathcal{X}\otimes _{\mathcal{F}}\Re=
\mathcal{X}\overline{\otimes }\Re$. By Lemma \ref{A25}(a), if $\mathcal{Y}$
and $\mathcal{Z}$ are two dual operator spaces and $\varphi $ is a
weak*-continuous completely bounded map from $\mathcal{Y}$ to $\mathcal{Z}$,
then there exists a weak*-continuous completely bounded map $\mathrm{id}
\otimes \varphi $ from $\mathcal{X}\otimes _{\mathcal{F}}\mathcal{Y}$ to $
\mathcal{X}\otimes _{\mathcal{F}}\mathcal{Z}$ such that $(\mathrm{id}\otimes
\varphi )(t)(\omega \otimes \nu )=t(\omega \otimes \varphi _{\ast }(\nu ))$ (
$t\in \mathcal{X}\otimes _{\mathcal{F}}\mathcal{Y}$; $\omega \in \mathcal{X}
_{\ast }$; $\nu \in \mathcal{Z}_{\ast }$). This enables us to define the
following.

\medskip
\begin{definition}
\label{C2} A weak*-continuous completely bounded map $\beta $ from $
\mathcal{X}$ to $\mathcal{X}\otimes _{\mathcal{F}}\Re$
(respectively, $\Re\otimes _{\mathcal{F}}\mathcal{X}$) is said to be
a \emph{normal right coaction} (respectively, {\em normal left coaction}) if $
(\beta \otimes \mathrm{id})\beta =(\mathrm{id}\otimes \delta
)\beta $ (respectively, $(\mathrm{id}\otimes \beta )\beta
=(\delta \otimes \mathrm{id})\beta $).
\end{definition}

\medskip
\begin{remark}
\label{C3} (a) We recall from \cite{Ruan2} that a \emph{right operator $
\Re_{\ast }$-module} is an operator space $N$ with a completely
bounded map $m$ from $N\hat{\otimes}\Re_{\ast }$ to $N$ such that $
m\circ (m\otimes \mathrm{id})=m\circ (\mathrm{id}\otimes \delta _{\ast })$
\emph{(left operator $\Re_{\ast }$-module} can be defined
similarly). For any normal right coaction $\beta $, the predual map $
\beta _{\ast }$ from $\mathcal{X}_{\ast }\hat{\otimes}\Re
_{\ast }$ to $\mathcal{X}_{\ast }$ gives a right operator $\Re_{\ast
}$-module structure on $\mathcal{X}_{\ast }$.

\smallskip
\noindent (b) It is natural to ask whether the dual $\Re_{\ast }$-module structure on $\mathcal{X}$ comes from a normal left coaction
(on $\mathcal{X}^{\ast }$). However, it is not clear why this
$\Re_*$-multiplication can be extended to the operator
projective tensor product (or
the range of the dual map lies in the Fubini product). Nevertheless, we will
see later that it comes from a more general form of coaction (Lemma \ref{C8}
).
\end{remark}

\medskip
\noindent
{\bf Notation:}
{\em Throughout this section, $\Re^{n}$ is the n-th times von Neumann
algebra tensor product of $\Re$ whereas $\Re_{\ast }^{n}$ is
the n-th times operator projective tensor product of $\Re_{\ast }$ ($
n\geq 1)$ and we take $\Re^{0}=\mathbb{C}=\Re_{\ast
}^{0}$.}

\bigskip
Suppose that $\mathcal{X}$ is a dual operator space with normal right coaction $
\beta $ and normal left coaction $\gamma $ such that $(\mathrm{id}
\otimes \beta )\circ \gamma =(\gamma \otimes \mathrm{id})\circ
\beta $. With the help of Lemma \ref{A25}(a), we can define as in the
case of Hopf $C^{\ast }$-algebra, a map $\delta _{n}$ from $\mathcal{X}
\otimes _{\mathcal{F}}\Re^{n}$ to $\mathcal{X}\otimes _{\mathcal{F}}
\Re^{n+1}$ by $\delta _{n}(x\otimes s_{1}\otimes ...\otimes
s_{n})=\beta (x)\otimes s_{1}\otimes ...\otimes
s_{n}+\sum_{k=1}^{n}(-1)^{k}x\otimes s_{1}\otimes ...\otimes s_{k-1}\otimes
\delta (s_{k})\otimes s_{n+1}\otimes ...\otimes s_{n}+(-1)^{n+1}(\gamma
(x)\otimes s_{1}\otimes ...\otimes s_{n})^{\sigma _{n+1,1}}$ ($n\geq 1$) and
$\delta _{0}(x)=\beta (x)-\gamma (x)^{\sigma _{1,1}}$ (where $\sigma
_{n,k}$ is the map from $\Re^{k}\otimes _{\mathcal{F}}\mathcal{X}
\otimes _{\mathcal{F}}\Re^{n-k}$ to $\mathcal{X}\otimes _{\mathcal{F}
}\Re^{n}$ as defined in Section 2). The same argument as Lemma
\ref {B8} shows that this gives a cochain complex and the cohomology
defined is called the \emph{normal natural cohomology }of
$\Re$ with coefficient in ($\mathcal{X},\beta ,\gamma )$.

\medskip

On the other hand, we can also define a map $\partial _{n}$ from $\mathrm{CB}
_{\sigma }(\mathcal{X};\Re^{n})$ (the set of all weak*-continuous
completely bounded maps from $\mathcal{X}$ to $\Re^{n}$) to $\mathrm{
CB}_{\sigma }(\mathcal{X};\Re^{n+1})$ by $\partial _{n}(T)=(T\otimes
\mathrm{id})\circ \beta +\sum_{k=1}^{n}(-1)^{k}(\mathrm{id}
^{n-k}\otimes \delta \otimes \mathrm{id}^{k-1})\circ T+(-1)^{n+1}(\mathrm{
id\otimes T)\circ \gamma }$ ($n\geq 1)$ and $\partial _{0}(f)=(f\otimes
\mathrm{id})\circ \beta -(\mathrm{id}\otimes f\mathrm{)\circ \gamma }$. 
As in the case of Hopf $C^*$-algebra, $(\mathrm{CB}_{\sigma }(
\mathcal{X};\Re^{n}),\partial _{n})$ is a cochain complex and
induces a cohomology $H_{\sigma ,d}^{n}(\Re;\mathcal{X})=\mathrm{Ker}
(\partial _{n})/\mathrm{Im}(\partial _{n-1})$ which is called the \emph{
normal dual cohomology of }$\Re$ with coefficient in $(\mathcal{X}
,\beta ,\gamma )$. In the case when $\gamma =0$, we denote it by $
H_{\sigma ,d,r}^{n}(\Re;\mathcal{X})$. Now using a similar argument
as that of Proposition \ref{B12}, we have the following.

\medskip
\begin{proposition}
\label{C4}
Let $\Re$ be a saturated Hopf von Neumann algebra as above.

\smallskip
\noindent (a) If $\Re_{\ast }$ is unital, then $H_{\sigma ,d,r}^{n}(
\Re;\mathcal{X})=(0)$ ($n=1,2,3,...$) for any dual operator space $
\mathcal{X}$ with a normal right coaction by $\Re$.

\smallskip
\noindent (b) If $H_{\sigma ,d,r}^{1}(\Re;\Re)=(0)$, then $
\Re_{\ast }$ is unital.
\end{proposition}

\medskip
However, we are more interested in the existence of a bounded approximate
identity of $\Re_{\ast }$ (which is related to amenability). A
closer look at the above reveals that this can be achieved if we remove the
weak*-continuity. Moreover, by doing so, we can also extend the definition of
coactions to general operator spaces.
Let us first note
that $\mathcal{X}\otimes _{\mathcal{F}}\Re=(\mathcal{X}_{\ast }\hat{
\otimes}\Re_{\ast })^{\ast }=\mathrm{CB}(\Re_{\ast };
\mathcal{X})$. Now by translating the coaction identity in terms of $
\mathrm{CB}(\Re_{\ast };\mathcal{X})$, we can define a more
general form of coactions.

\medskip
\begin{definition}
\label{C5} A completely bounded map $\beta $ from $X$ to $\mathrm{CB}(
\Re_{\ast };X)$ is said to be a \emph{right coaction} (respectively, a
\emph{left coaction}) if $\beta (\beta (x)(\omega ))(\nu
)=\beta (x)(\delta _{\ast }(\nu \otimes \omega ))$ (respectively, $
\beta (\beta (x)(\omega ))(\nu )=\beta (x)(\delta _{\ast
}(\omega \otimes \nu ))$) for any $x\in X$ and $\omega ,\nu \in \Re
_{\ast }$. Moreover, a right (or a left) coaction $\beta $ is said
to be \emph{non-degenerate} if $\beta (X)(\Re_{\ast })$ is
dense in $X$ (which is equivalent to weakly dense).
We call $(X,\beta ,\gamma )$ a
{\em $\Re$-bicomodule} if $\beta $ is a right coaction and $\gamma $
is a left coaction on $X$ by $\Re$ such that $\beta (\gamma
(x)(\omega
))(\nu )=\gamma (\beta (x)(\nu ))(\omega )$ for any $x\in X$ and $
\omega ,\nu \in \Re_{\ast }$.
\end{definition}

\medskip
Note that $\beta (\beta (x)(\omega ))(\nu )=\beta
^{\#}(\beta )(x)(\omega )(\nu )$ and $\beta (x)(\delta _{\ast
}(\omega \otimes \nu ))=\beta ^{0}(\delta _{\ast })(x)(\omega \otimes
\nu )$ (where $\beta ^{\#}$ and $\beta ^{0}$ are the maps as
defined in Lemma \ref{A27}(b)). Hence the right (respectively,
left) coaction identity in Definition \ref{C5} can be simplified to $
\beta ^{\#}(\beta )=\beta ^{0}(\delta _{\ast }\circ \sigma
)$ (respectively, $\beta ^{\#}(\beta )=\beta ^{0}(\delta
_{\ast })$) under the standard identification $\mathrm{CB}(\Re_{\ast
};\mathrm{CB}(\Re_{\ast };X))\cong \mathrm{CB}(\Re_{\ast }
\hat{\otimes}\Re_{\ast };X)$ (see Remark \ref{A26}$)$. Note that the
forms of these simplified coaction identities depend on whether we take the
standard identification or the reverse identification (see Remark \ref{A26}).

\medskip
\begin{example}
\label{C19} (a) Suppose that $\beta $ is a coaction of $\Re$
on a von Neumann algebra $\mathcal{M}$. Let $N$ be any subset of $\mathcal{M}
$ and $X_{N}$ be the closed linear span of the set $\{(\mathrm{id}\otimes
\omega )(\beta (x)):\omega \in \Re_{\ast };x\in N\}$. We
denote by $\beta _{N}$ the composition of the restriction of $
\beta $ on $X_{N}$ with the complete isometry from $\mathcal{M}\bar{
\otimes}\Re$ to $\mathrm{CB}(\Re_{\ast };\mathcal{M})$. Then
we have $\beta _{N}[(\mathrm{id}\otimes \omega )(\beta (x))](\nu
)=(\mathrm{id}\otimes (\nu \cdot \omega ))(\beta (x))\in X_{N}$ (for
any $\omega ,\nu \in \Re_{\ast }$ and $x\in N$) and it is not hard
to see that $\beta _{N}$ is a right coaction on $X_{N}$.

\smallskip
\noindent (b) Suppose that $\Re$ comes from a Kac algebra $\mathbf{K}
$ and $\beta $ is any completely contractive right coaction of $\Re$
on any operator space $X$. Let $N$ and $U$ be the unit balls of $X$ and $
\Re_{\ast }$ respectively. Since $\beta \in \mathrm{CB}(X;
\mathrm{CB}(\Re_{\ast };X))\cong \mathrm{CB}(\Re_{\ast };
\mathrm{CB}(X;X))$ is a complete contraction, $\Vert \beta (\omega
)(x)\Vert \leq 1$ for any $\omega \in U$ and $x\in N$. It is not hard to see
that this defines an action of $\mathbf{K}$ on $N$ in the sense of \cite[2.2]
{ES}.

\smallskip
\noindent (c) For any Hilbert space $H$, there is an one to one
correspondence between right coactions of $\Re$ on the column Hilbert
space $H_{c}$ and the representations of $\Re$ on $H$ (see
\cite[2.11]{Ng3}).

\smallskip
\noindent (d) Let $\Gamma$ be a discrete group and $\beta$ be a coaction
of $C^*_r(\Gamma)$ on a $C^*$-algebra $A$.
Then $\beta$ induces a right
coaction $\bar\beta$ of the group von Neumann algebra $vN(\Gamma)$ on $A$
(by $\bar\beta(a)(\omega) = (\id\otimes \omega)\beta(a)$ for any $a\in A$
and $\omega\in A(\Gamma) = vN(\Gamma)_*$).
Thus, by the next Lemma, we have a left
coaction $\check\beta$ of $vN(\Gamma)$ on $A^*$.
Moreover, as $\beta$ is
injective and $A^*\otimes A(\Gamma)$ separates points of $M(A\otimes
C^*_r(\Gamma))$, the subspace
$\check\beta(A^*)(A(\Gamma))$ is weak*-dense in $A^*$.
Notice that the function $\varphi_s$ defined in Example \ref{B1}(b) is in
$A(\Gamma)$.
For any $t\in
\Gamma$, the sets $\{ f\in A^*:
\check\beta(f)(\varphi_s) = \delta_{s,t} f\}$ (where $\delta_{s,t}$ is the
Kronecker delta) and $\{ \check\beta(g)(\varphi_t): g\in A^*\}$ coincide
and we denote this set by $A^*_t$.
If $r\neq s\in \Gamma$, $a\in A_r$ and
$f\in A^*_s$, then $f(a) = \check\beta(f)(\varphi_s)(a) = f(\id\otimes
\varphi_s)
\beta(a)=0$.

\end{example}

\medskip
\begin{lemma}
\label{C8} A right coaction $\beta $ of $\Re$
on $X$ induces a left coaction $\check{
\beta}$ on $X^{\ast }$ such that $\check{\beta}(f)(\omega
)(x)=f(\beta (x)(\omega ))$ ($f\in X^{\ast }$; $x\in X$; $\omega \in
\Re_{\ast }$). Similarly, a left coaction on $X$ will induce a
right coaction on $X^{\ast }$.
\end{lemma}
\noindent \textbf{Proof:} Let $\check{\beta}$ be the composition of
the completely bounded map $\beta ^{\#}:X^{\ast }\longrightarrow
\mathrm{CB}(X;\Re)$ (see Lemma \ref{A27}(b)) with the canonical
complete isometry from $\mathrm{CB}(X;\Re)$ to $\mathrm{CB}(
\Re_{\ast };X^{\ast })$ (see Lemma \ref{A25}(d)). Thus, $\check{\beta}$
is completely bounded and $\check{\beta}(f)(\omega )(x)=\beta
^{\#}(f)(x)(\omega )=f(\beta (x)(\omega ))$. It remains to show the
left coaction identity.
Indeed, $\check{\beta}(\check{\beta}(f)(\omega
))(\nu )(x)=f(\beta (\beta (x)(\nu ))(\omega ))=f(\beta
(x)(\delta _{\ast }(\omega \otimes \nu )))=\check{\beta}(f)(\delta
_{\ast }(\omega \otimes \nu ))(x)$
for any $f\in X^{\ast }$, $x\in X$ and $\omega ,\nu
\in \Re_{\ast }$. The proof of the second statement is the
same.

\bigskip

In fact, by a similar argument as in Example \ref{C19}(d), the left
coaction $\check\beta$ is normal. 
It is ``weakly non-degenerate'' if
$\beta$ is injective.
However, we will not need these facts in this paper

\medskip

We will again define two cohomology theories for this type of bicomodules.
We first consider the analogue of the natural cohomology.
Suppose that $\beta $ is a normal right coaction of $\Re$
on the dual operator space $\mathcal{X}$.
As $\mathrm{CB}(\Re_{\ast }^{n};
\mathcal{X})\cong \mathcal{X}\otimes _{\mathcal{F}}\Re^{n}$ under
the identification: $T_{\alpha }(\omega _{1}\otimes ...\otimes \omega _{n})=(
\mathrm{id}\otimes \omega _{1}\otimes ...\otimes \omega _{n})(\alpha )$
(for $\alpha \in \mathcal{X}\otimes _{\mathcal{F}}\Re^{n}$), the map
$\beta _{(n)}\in \mathrm{CB}(\mathcal{X}\otimes _{\mathcal{F}}\Re^{n};
\mathcal{X}\otimes _{\mathcal{F}}\Re^{n+1})$ defined by $
\beta _{(n)}(z)=(\beta \otimes \mathrm{id}^{n})(z)$ ($z\in
\mathcal{X}\otimes_{\cal F} \Re^{n}$) can be identified with the map from $
\mathrm{CB}(\Re_{\ast }^{n};\mathcal{X})$ to $\mathrm{CB}(\Re
_{\ast }^{n+1};\mathcal{X})$ given by $\beta _{(n)}(T)(\omega
_{1}\otimes ....\otimes \omega _{n+1})=\beta (T(\omega _{2}\otimes
....\otimes \omega _{n+1}))(\omega _{1})$ ($T\in \mathrm{CB}(\Re
_{\ast }^{n};\mathcal{X})$).

\medskip

Now for a general $\Re$-bicomodule $(X,\beta ,\gamma )$, let $
\beta _{(n)}$, $\gamma _{(n)}$ and $\delta _{n,k}$ be maps from $
\mathrm{CB}(\Re_{\ast }^{n};X)$ to $\mathrm{CB}(\Re_{\ast
}^{n+1};X)$ given by $\beta _{(n)}(T)(\omega _{1}\otimes ....\otimes
\omega _{n+1})=\beta (T(\omega _{2}\otimes ....\otimes \omega
_{n+1}))(\omega _{1})$, $\gamma _{(n)}(T)(\omega _{1}\otimes ....\otimes
\omega _{n+1})=\gamma (T(\omega _{1}\otimes ....\otimes \omega _{n}))(\omega
_{n+1})$ and $\delta _{n,k}(T)=T\circ (\mathrm{id}^{k-1}\otimes \delta
_{\ast }\otimes \mathrm{id}^{n-k})$. Then $\beta _{(n)}$ is completely
bounded since it is the map $\tilde{\beta}$ in Lemma \ref{A27}(a)
under the standard identification of Remark \ref{A26}. 
The same is true for $\gamma _{(n)}$.
Let
$$\delta _{n}=
\cases{
\beta_{(n)}+\sum_{k=1}^n(-1)^k\delta _{n,k}+(-1)^{n+1}\gamma _{(n)} &
\quad $n\geq 1$\cr
\beta (x)-\gamma (x) & \quad $n=0$.
}
$$
A proof is needed to show that $\delta_n$ gives a cochain complex.

\medskip
\begin{lemma}
\label{C11} $\delta _{n}\circ \delta _{n-1}=0$ ($n=1,2,3...$).
\end{lemma}
\noindent \textbf{Proof:} For $n=1$, we have $\delta _{1}(\delta
_{0}(x))(\omega \otimes \nu )=\beta \lbrack \beta (x)(\nu
)-\gamma (x)(\nu )](\omega )-\beta (x)(\delta _{\ast }(\omega \otimes
\nu ))+\gamma (x)(\delta _{\ast }(\omega \otimes \nu ))+\gamma \lbrack
\beta (x)(\omega )-\gamma (x)(\omega )](\nu )=0$ ($x\in X;$ $\omega
,\nu \in \Re_{\ast }$) (by the left and the right coaction
identities).
For $n>1$,
\[
\beta _{(n)}\circ \beta _{(n-1)}=\delta _{n,1}\circ \beta
_{(n-1)},\quad \beta _{(n)}\circ \delta _{n-1,k}=\delta _{n,k+1}\circ
\beta _{(n-1)},\quad \beta _{(n)}\circ \gamma _{(n-1)}=\gamma
_{(n)}\circ \beta _{(n-1)},
\]
\[
\gamma _{(n)}\circ \delta _{n-1,k}=\delta _{n,k}\circ \gamma _{(n-1)}\quad
\mathrm{and}\quad \gamma _{(n)}\circ \gamma _{(n-1)}=\delta _{n,n}\circ
\gamma _{(n-1)}.
\]
Thus, in order to show $\delta _{n}\circ \delta _{n-1}=0$, we need to check
that $\sum_{l=1}^{n}\sum_{k=1}^{n-1}(-1)^{l+k}\delta _{n,l}\circ \delta
_{n-1,k}=0$. This can be shown again by a decomposition and a comparison
(similar to Lemma \ref{B8}).

\bigskip

As in section 2, we call the cohomology $H^{n}(\Re;X)=\mathrm{Ker}
(\delta _{n})/\mathrm{Im}(\delta _{n-1})$ ($n=1,2,3,...)$ and $H^{0}(
\Re;X)=\{x\in X:\beta (x)=\gamma (x)\}$ the \emph{natural
cohomology }\emph{of }$S$ \emph{with coefficient in }$(X,\beta ,\gamma
)$.

\medskip

Next, we want to define the dual cohomology analogue for Hopf von
Neumann algebras.
By Lemma \ref{A27}(b), $\beta $ induces a completely bounded map $\beta _{n}$
given by $\beta _{n}(F)=\beta ^{\#}(F)\in \mathrm{CB}(X;\mathrm{CB}(
\Re_{\ast };\Re^{n}))=\mathrm{CB}(X;\Re^{n}\bar{
\otimes}\Re)$
(for any $F\in \mathrm{CB}(X;\Re^{n})$),
i.e. $\beta _{n}(F)(x)(\omega _{1}\otimes
...\otimes \omega _{n+1})=F(\beta (x)(\omega _{n+1}))(\omega
_{1}\otimes ...\otimes \omega _{n})$ ($x\in X$; $\omega _{1},...,\omega
_{n+1}\in \Re_{\ast }$). Similarly, $\gamma $ induces a completely
bounded map $\gamma _{n}$ such that $\gamma _{n}(F)(x)(\omega _{1}\otimes
...\otimes \omega _{n+1})=F(\gamma (x)(\omega _{1}))(\omega _{2}\otimes
...\otimes \omega _{n+1})$. Now we let $\partial _{n,k}(F)=(\mathrm{id}
^{n-k}\otimes \delta \otimes \mathrm{id}^{k-1})\circ F$ and
$$\partial_{n}=
\cases{
\beta _{n}+\sum_{k=1}^{n}(-1)^{k}\partial _{n,k}+(-1)^{n+1}\gamma
_{n} & \quad $n\geq 1$\cr
\beta _{0}-\gamma _{0} & \quad $n=0$.
}
$$
By Lemma
\ref{C11} and the proof of Proposition \ref{C10} below, we have the
following.

\medskip
\begin{lemma}
\label{C12} $\partial _{n}\circ \partial _{n-1}=0$ for $n=1,2,...$.
\end{lemma}

\medskip
Thus, $\{\partial _{n}\}$ defines a cohomology $H_{d}^{n}(\Re;X)=
\mathrm{Ker}(\partial _{n})/\mathrm{Im}(\partial _{n-1})$ ($n=1,2,3,...)$
and $H_{d}^{0}(\Re;X)=\{f\in X^{\ast }:\beta ^{\#}(f)=\gamma
^{\#}(f)\}$ which is called the \emph{dual cohomology }\emph{of }$S$ \emph{
with coefficient in }$X$. The use of the term ``dual cohomology'' can be
justified by the following equivalent formulation. This also shows that $
H^{n}(\Re,\bullet )$ is a more general form of cohomology theory than $
H_{d}^{n}(\Re,\bullet )$.

\medskip

Suppose that $(X,\beta ,\gamma )$ is a $\Re$-bicomodule and $
\check{\beta}$ and $\check{\gamma}$ are respectively the left and the
right coactions on $X^{\ast }$ given by Lemma \ref{C8}. Then it is easy to
see that $(X^{\ast },\check{\beta},\check{\gamma})$ is again a $
\Re$-bicomodule.

\medskip
\begin{proposition}
\label{C10}
For any saturated Hopf von Neumann algebra $\Re$ and any
$\Re$-bicomodule $X$, we have $H_{d}^{n}(\Re;X)\cong
H^{n}(\Re;X^{\ast })$ (for $n=0,1,2,3,...$).
\end{proposition}
\noindent \textbf{Proof:} The idea of proof rely on the fact that $\mathrm{CB
}(\Re_{\ast }^{n};X^{\ast })\cong \mathrm{CB}(X;\Re^{n})$.
In this case, the corresponding element $\overline{T}$ of $T\in \mathrm{CB}(
\Re_{\ast }^{n};X^{\ast })$ is given by $\overline{T}(x)(\omega
_{1}\otimes ...\otimes \omega _{n})=T(\omega _{1}\otimes ...\otimes \omega
_{n})(x)$. Thus,
\begin{eqnarray*}
\overline{\check{\beta}_{(n)}(T)}(x)(\omega _{1}\otimes ...\otimes
\omega _{n+1}) &=&\check{\beta}(T(\omega _{1}\otimes ...\otimes \omega
_{n}))(\omega _{n+1})(x)\quad =\quad T(\omega _{1}\otimes ...\otimes \omega
_{n})(\beta (x)(\omega _{n+1})) \\
&=&\beta _{n}(\overline{T})(x)(\omega _{1}\otimes ...\otimes \omega
_{n+1})
\end{eqnarray*}
\noindent
(note that $\check\beta$ is a left coaction).
Similarly, $\overline{\check{\gamma}_{(n)}(T)}=\gamma _{n}(\overline{T})$.
On the other hand, it is easy to see that $\overline{\delta _{n,k}(T)}
=\partial _{n,n-k+1}(\overline{T})$. Therefore, $\overline{\delta _{n}(T)}
=\gamma _{n}(\overline{T})+\sum_{l=1}^{n}(-1)^{n-l+1}\partial _{n,l}(
\overline{T})+(-1)^{n+1}\beta _{n}(\overline{T})=(-1)^{n+1}\partial
_{n}(\overline{T})$.

\bigskip

Next, we would like to study the vanishing of the dual cohomology of Hopf
von Neumann algebras. Again, let us first consider the one-sided case when $
\gamma =0$ and use $H_{d,r}^{n}$ to denote the dual cohomology defined in
this situation.

\medskip
\begin{theorem}
\label{C6} Let $\beta $ be any right coaction of the saturated Hopf
von Neumann algebra $\Re$ on an operator space $X$.

\smallskip
\noindent (a) $H_{d,r}^{0}(\Re;X)=(0)$ if and only if $\beta $
is non-degenerate.

\smallskip
\noindent (b) If $\Re_{\ast }$ has a bounded left approximate
identity, then $H_{d,r}^{n}(\Re;X)=(0)$ ($n\geq 1$).

\smallskip
\noindent (c) If $H_{d,r}^{1}(\Re;\Re)=(0)$, then
$\Re_{\ast }$ has a bounded left approximate identity.
\end{theorem}
\noindent \textbf{Proof:} (a) This part is clear.

\smallskip
\noindent (b) Suppose that $\{\nu _{i}\}$ is a bounded left approximate
identity of $\Re_{\ast }$ and $T\in \mathrm{Ker}(\partial _{n})$.
Consider the following identification: $\mathrm{CB}(X;\Re\bar{\otimes
}\Re^{n-1})\cong \mathrm{CB}(X;\mathrm{CB}(\Re_{\ast };
\Re^{n-1}))\cong \mathrm{CB}(\Re_{\ast };\mathrm{CB}(X;
\Re^{n-1}))$ (note that the first isomorphism is different from the
one considered in the paragraph preceding Lemma \ref{C12}). Let $\hat{T}\in
\mathrm{CB}(\Re_{\ast };\mathrm{CB}(X;\Re^{n-1}))$ be the
corresponding element of $T$ (i.e. $\hat{T}(\omega _{0})(x)(\omega
_{1}\otimes ...\otimes \omega _{n-1})=T(x)(\omega _{0}\otimes ...\otimes
\omega _{n-1})$). Since $\partial _{n}(T)=0$, we have, for any $\omega
_{0},...,\omega _{n}\in \Re_{\ast }$,
\begin{eqnarray*}
0 &=&\hat{T}(\omega _{0})(\beta (x)(\omega _{n}))(\omega _{1}\otimes
...\otimes \omega _{n-1})+\,\sum_{k=1}^{n-1}(-1)^{k}\hat{T}(\omega
_{0})(x)(\omega _{1}\otimes ...\otimes (\omega _{n-k}\cdot \omega
_{n-k-1})\otimes ...\otimes \omega _{n})\,+ \\
&&\quad \quad \quad (-1)^{n}\hat{T}(\omega _{0}\cdot \omega _{1})(x)(\omega
_{2}\otimes ...\otimes \omega _{n}).
\end{eqnarray*}
\noindent Moreover, as $\mathrm{CB}(X;\Re^{n-1})\cong (X\hat{\otimes}
\Re_{\ast }^{n-1})^{\ast }$ (see Lemma \ref{A25}(c)),
the bounded
net $\{\hat{T}(\nu _{i})\}$ has a subnet $\{\hat{T}(\nu _{i_{j}})\}$ that
weak*-converges to some $F\in \mathrm{CB}(X;\Re^{n-1})$. Note that $
\nu _{i_{j}}\cdot \omega $ converges to $\omega $ for any $\omega \in
\Re_{\ast }$. Therefore, by putting $\omega _{0}=\nu _{i_{j}}$ into
the above equation and taking limit, we obtain $0=F(\beta (x)(\omega
_{n}))(\omega _{1}\otimes ...\otimes \omega
_{n-1})+\sum_{k=1}^{n-1}(-1)^{k}F(x)(\omega _{1}\otimes ...\otimes (\omega
_{n-k}\cdot \omega _{n-k-1})\otimes ...\otimes \omega _{n})+(-1)^{n}\hat{T}
(\omega _{1})(x)(\omega _{2}\otimes ...\otimes \omega _{n})$ and so $
T=(-1)^{n-1}\partial _{n-1}(F)$ as required.

\smallskip
\noindent (c) Recall that the right coaction $\beta $ of $\Re$ on $
\Re$ is given by $\beta (s)(\omega )=(\mathrm{id}\otimes
\omega )\delta (s)$ ($s\in \Re;$ $\omega \in \Re_{\ast }$).
As in the proof of Proposition \ref{B12}(b), because $\mathrm{id}\in \mathrm{
Ker}(\partial _{1})$, there exists $u\in \Re^{\ast }$ such that $
\partial _{0}(u)=\mathrm{id}$. Thus, $(u\times _{2}\omega )(s)=u((\mathrm{id}
\otimes \omega )\delta (s))=u(\beta (s)(\omega ))=\partial
_{0}(u)(s)(\omega )=s(\omega )$ for any $\omega \in \Re_{\ast }$ and
$s\in \Re$ (where $\times _{2}$ is the second Arens product on $
\Re^{\ast }=(\Re_{\ast })^{\ast \ast }$). Thus $\Re
^{\ast }$ has a left identity for the second Arens product and so
$\Re_{\ast }$ has a bounded left approximate identity (see e.g.
\cite[5.1.8]{Pal}).

\bigskip

By Proposition \ref{C10}, $H_{d,r}^{n}(\Re;X)\cong H_{l}^{n}(
\Re;X^{\ast })$ (where $H_{l}^{n}$ is the natural cohomology
obtained in the case when the right coaction is zero) for any
right $\Re$-comodule $X$. The following corollary shows that the vanishing of $H_{l}^{n}
$ lies between the existence of a bounded left approximate identity and the
existence of an identity in $\Re_{\ast }$. In the case when $
\Re$ is the bidual of a Hopf $C^{\ast }$-algebra, all these three
properties coincide (by \cite[2.6]{Ng2}).

\medskip
\begin{corollary}
\label{C13} If $H_{l}^{1}(\Re;\Re^{\ast })=(0)$, then there
exists a bounded left approximate identity for $\Re_{\ast }$.
On the other hand, if $\Re_{\ast }$ is unital, then
$H_{l}^{n}(\Re;X)=(0)$ for any $n\geq 1$ and any left
$\Re$-comodule $X$.
\end{corollary}

\medskip
In fact, the first statement follows clearly from Proposition \ref{C10} and
Theorem \ref{C6}(c). Moreover, suppose that $u$ is the identity of
$\Re_{\ast }$ and $T\in \mathrm{Ker}(\delta _{n+1})$.
If $F\in \mathrm{CB}(
\Re_{\ast }^{n-1};X^{\ast })$ is defined by $F(\omega _{1}\otimes
...\otimes \omega _{n-1})=T(u\otimes \omega _{1}\otimes ...\otimes \omega
_{n-1})$, then it is not hard to see that $T=\delta _{n}((-1)^{n}F)$.

\medskip

It turns out that there is no need to study the vanishing of the 2-sided
dual cohomology of a Hopf von Neumann algebra because of its relation with
operator cohomology (that studied in \cite{Ruan2}). Let us first recall
from \cite{Ruan2} the definition of operator cohomology. Suppose that $B$ is
a completely contractive Banach algebra and $V$ is an operator $B$-bimodule
(see \cite[p.1453]{Ruan2} or Remark \ref{C3}(a)).
Consider a map $d_{n}$ from $B^{n}\hat{\otimes}V$
to $B^{n-1}\hat{\otimes}V$ (where $B^{n}$ is the n-th times operator
projective tensor product of $B$ and $B^{0}=\mathbb{C}$) given by
\begin{eqnarray*}
\lefteqn{d_{n}(a_{1}\otimes ...\otimes a_{n}\otimes v)=}\\
& & \quad a_{1}\otimes ...\otimes
a_{n-1}\otimes (a_{n}\cdot v)+\sum_{i=1}^{n-1}(-1)^{n-i}a_{1}\otimes
...\otimes a_{i}\cdot a_{i+1}\otimes ...\otimes a_{n}\otimes
v+(-1)^{n}a_{2}\otimes ...\otimes a_{n}\otimes v\cdot a_{1}.
\end{eqnarray*}
Note that the $d_{n}$'s that we use here differ from those
in \cite{Ruan2} by a sign of $(-1)^{n}$.
Now $((B^{n}\hat{\otimes}V)^{\ast },d_{n}^{\ast })$ is a cochain
complex and defines the \emph{operator cohomology,} $\mathrm{OH}
^{n}(B;V^{\ast })$, from $B$ to $V^{\ast }$ (i.e. $\mathrm{
OH}^{n}(B;V^{\ast })=\mathrm{Ker(}d_{n+1}^{\ast })/\mathrm{Im(}d_{n}^{\ast
}))$ (see \cite[p.1455]{Ruan2}).

\medskip

Consider now the completely contractive Banach algebra $\Re_{\ast }
$. Since $\mathrm{CB}(\Re_{\ast }\hat{\otimes}X;X)\cong \mathrm{CB}
(X;\mathrm{CB}(\Re_{\ast };X))$, if $m\in \mathrm{CB}(\Re
_{\ast }\hat{\otimes}X;X)$ is a completely bounded left $\Re_{\ast }$-multiplication on $X$, then the corresponding map $\beta _{m}\in
\mathrm{CB}(X;\mathrm{CB}(\Re_{\ast };X))$ given by $\beta
_{m}(x)(\omega )=m(\omega \otimes x)$ ($x\in X;$ $\omega \in \Re$)
is a right coaction.
Using the same formula, a right coaction defines a left $\Re_{\ast }$-multiplication and these give the following.

\medskip
\begin{lemma}
\label{C22}There is an one to one correspondence between completely bounded
left (respectively, right) $\Re_{\ast }$-module structures on $X$
and right (respectively, left) coactions of $\Re$ on $X$.
\end{lemma}

\medskip
This implies that there is an one to one correspondence between operator $
\Re_{\ast }$-bimodule structures and $\Re$-bicomodule
structures on $X$.
Recall that if $X$ is a $\Re$-bicomodule, the $
\Re_{\ast }$-multiplication on $X$ is given by $\omega \cdot
x=\beta (x)(\omega )$ and $x\cdot \omega =\gamma (x)(\omega )$ ($
\omega \in \Re_{\ast };$ $x\in X)$. We can now show that the cochain
complex that defines $H_{d}^{n}(\Re;X)$ is the same as the one that
defines the operator cohomology from $\Re_{\ast }$ to $X^*$.

\medskip
\begin{proposition}
\label{C15} For any saturated Hopf von Neumann algebra $\mathcal{R}$ and
any $\Re$-bicomodule $X$, we have $H_{d}^{n}(\Re;X)=
\mathrm{OH}^{n}(\Re_{\ast };X^{\ast })$ (for $n=1,2,3...$).
\end{proposition}
\noindent \textbf{Proof: }Note that $\mathrm{CB}(X;\Re^{n})\cong (
\Re_{\ast }^{n}\hat{\otimes}X)^{\ast }$ under the identification $
f_{T}(\omega _{1}\otimes ...\otimes \omega _{n}\otimes x)=T(x)(\omega
_{1}\otimes ...\otimes \omega _{n})$ ($T\in \mathrm{CB}(X;\Re^{n})$
). It is clear that for any $\omega _{1},...,\omega _{n+1}\in $ $\Re
_{\ast }$ and $x\in X$,
\[
\beta _{n}(T)(x)(\omega _{1}\otimes ...\otimes \omega
_{n+1})=f_{T}(\omega _{1}\otimes ...\otimes \omega _{n}\otimes (\omega
_{n+1}\cdot x)),
\]
\[
\partial _{n,k}(T)(x)(\omega _{1}\otimes ...\otimes \omega
_{n+1})=f_{T}(\omega _{1}\otimes ...\otimes \omega _{n-k}\otimes \omega
_{n-k+1}\cdot \omega _{n-k+2}\otimes \omega _{n-k+3}\otimes ...\otimes
\omega _{n+1}\otimes x)
\]
and
\[
\gamma _{n}(T)(x)(\omega _{1}\otimes ...\otimes \omega _{n+1})=f_{T}(\omega
_{2}\otimes ...\otimes \omega _{n+1}\otimes (x\cdot \omega _{1})).
\]
Hence $\partial _{n}(T)(x)(\omega _{1}\otimes ...\otimes \omega
_{n+1})=f_{T}(d_{n+1}(\omega _{1}\otimes ...\otimes \omega _{n+1}\otimes x))$
and the complexes $((\Re_{\ast }^{n}\hat{\otimes}X)^{\ast
},d_{n+1}^{\ast })$ and $(\mathrm{CB}(X;\Re^{n}),\partial _{n})$
coincide under the above identification.

\bigskip

This, together with \cite[2.1]{Ruan2}, gives the following
characterisation of the vanishing of the dual cohomology of Hopf von Neumann
algebras.

\medskip
\begin{corollary}
$H_{d}^{n}(\Re;X)=0$ for any $\Re$-bicomodule $X$ and any $
n\in \mathbb{N}$ if and only if $\Re_{\ast }$ is operator amenable.
\end{corollary}

\bigskip

\bigskip

\section{Dual cohomology and amenability}

\bigskip

In this section, we will give some interesting consequences of the results
in the previous sections.
If $A$ is a non-zero 2-sided $S$-invariant closed subalgebra of the dual
space $S^{\ast }$ (recall that $f\cdot g = (f\otimes g)\delta$ for
$f,g\in S^*$), then by \cite[III.2.7]{Take}, there exists a central
projection $e\in S^{\ast \ast }$ such that $A=eS^{\ast }$ and by
\cite[1.10(b)]{Ng1}, $A^{\ast }$ is a Hopf von Neumann algebra.

\bigskip
\noindent
{\bf Notation:}
{\em Throughout this section, we will assume that
$A$ is a non-zero 2-sided $S$-invariant closed subalgebra of
$S^{\ast }$ (note the difference between the $A$ in here and in
Section 1).
Moreover, $(R,\delta)$ is a not necessarily non-degenerate Hopf $C^*$-algebra.}

\medskip
\begin{lemma}
\label{D1} There is a complete contraction $\Psi _{A}$ from $M_{S}(X\otimes
S)$ to $\mathrm{CB}(A;X)$ such that $\Psi _{A}(m)(\omega )=(\mathrm{id}
\otimes \omega )(m)$ ($m\in M_{S}(X\otimes S);\omega \in A$). Moreover, if $A
$ separates points of $S$, then $\Psi _{A}$ is a complete isometry.
\end{lemma}
\noindent \textbf{Proof:} Let $j$ be the canonical $\ast $-homomorphism from
$S$ to $A^{\ast }$ (given by $j(s)(\omega)=\omega(s)$ for $s\in S$ and
$\omega\in A$).
Then Lemma \ref{1.12}(b) gives a complete contraction $
\mathrm{id}\otimes j$ from $M_{S}(X\otimes S)$ to $M_{j(S)}(X\otimes j(S))$.
Now by Proposition \ref{1.d} and the definition of the Fubini product, we
see that $M_{j(S)}(X\otimes j(S))$ can be regarded as an operator subspace of $
X^{\ast \ast }\otimes _{\mathcal{F}}A^{\ast }\cong \mathrm{CB}(A;X^{\ast
\ast })$ (note that for any $m\in M_{j(S)}(X\otimes j(S))$, $f\in X^*$
and $\omega\in A$, we have $(\id\otimes\omega)(m)\in X$ and
$(f\otimes\id)(m)\in M(j(S))$ by Lemma \ref{1.12}(a)).
It is not hard to see that the composition $\Psi _{A}$, of the
above two maps satisfies the required conditions
(in particular, $\Psi_A(M_S(X\otimes S))\subseteq \mathrm{CB}(A;X)$).
Furthermore, if $A$
separates points of $S$, then $j$ is a complete isometry
and so is $\Psi _{A}$.

\bigskip

If $\beta $ is a right coaction of $S$ on $X$, then the above lemma shows
that $\beta $ induces a right coaction $\bar{\beta}_{A}$ of $A^{\ast }$
on $X$ such that $\bar{\beta}_{A}(x)(\omega )=(\mathrm{id}\otimes
\omega )(\beta (x))$ for any $\omega \in A$ and $x\in X$ (the coaction
identity can be verified easily). Similarly, we can define from a left
coaction $\gamma $ of $S$ on $X$, a left coaction $\bar{\gamma}_{A}$
of $A^{\ast }$ on $X$.

\medskip
\begin{proposition}
\label{D8} Suppose that the dual space $S^{\ast }$ of a saturated
Hopf $C^*$-algebra $S$ contains a non-zero $S$-invariant closed
subalgebra $A$ that separates points of $S$. If $A$ is operator
amenable
(in the sense of \cite[2.2]{Ruan2}), then $H_{d}^{1}(S;X)=(0)$
for any $S$-bicomodule $X$.
Consequently, if $S^{\ast }$ is operator amenable, then
any first dual cohomology of $S$ vanishes.
\end{proposition}
\noindent
\textbf{Proof:} Let $(X,\beta ,\gamma )$ be a $S$-bicomodule. Suppose that $\bar{\gamma}_{A}$ and $\bar{\beta}_{A}$ are
the left and the right coactions of $A^{\ast }$ on $X$ as given above. By
the definition of operator amenability and Proposition \ref{C15}, $
H_{d}^{1}(A^{\ast };X)=(0)$. Consider $F\in \mathrm{CB}(X;M(S))$ such that $
\partial _{1}(F)=0$ and let $T=j\circ F$ $\in \mathrm{CB}(X;A^{\ast })$
(where $j$ is the map as in the proof of Lemma \ref{D1}). Note that $j$
``preserves the coproducts'' (i.e. $\delta(j(s)) = (j\otimes j)\delta(s)$
for all $s\in M(S)$).
For any $\omega ,\nu \in A$ and $x\in X$, by Lemma \ref{1.12}(c),
\begin{eqnarray*}
\lefteqn{T(\bar{\beta}_{A}(x)(\nu ))(\omega )-\delta (T(x))(\omega
\otimes \nu )+T(\bar{\gamma}_{A}(x)(\omega ))\left( \nu \right) } \\
&&\qquad \qquad \qquad =\quad \omega \circ F((\mathrm{id}\otimes \nu
)\beta (x))-(\omega \otimes \nu )\delta (F(x))+\nu \circ F((\omega
\otimes \mathrm{id})\gamma (x)) \\
&&\qquad \qquad \qquad =\quad (\omega \otimes \nu )\circ ((F\otimes \mathrm{
id})\circ \beta -\delta \circ F+(\mathrm{id}\otimes F)\circ \gamma
)(x)\quad =\quad 0.
\end{eqnarray*}
\noindent
Hence $T\in \mathrm{Ker}(\partial _{1})$ and there exists $f\in X^{\ast }$
such that for any $\omega \in A$ and $x\in X$, we have the
identities $\omega (F(x))=T(x)(\omega
)=f(\bar{\beta}_{A}(x)(\omega ))-f(\bar{\gamma}_{A}(x)(\omega
))=\omega ((f\otimes \mathrm{id})\beta (x)-(\mathrm{id}\otimes
f)\gamma (x))$. Since $A$ separates points of $S$ 
(and hence separates points of $M(S)$ as $A$ is $S$-invariant),
the above implies that $F=\partial _{0}(f)$.
Therefore, $H_{d}^{1}(S;X)=(0)$.

\bigskip

The proof of the above proposition actually shows that $H_{d}^{1}(S;X)
\subseteq H_{d}^{1}(A^{\ast };X)$ (when $A$ is a $S$-invariant closed
subalgebra of $S^{\ast }$ that separates points of $S$). Now, we would like
to consider a cohomology theory that is even ``smaller than'' $H_{d}^{1}(S;X)
$ (in fact, a ``restriction'' of it) that will help us to study the
vanishing of the dual cohomologies of the Hopf $C^{\ast }$-algebras
associated with locally compact groups. Let $U^{n}(S)=U_{\mathrm{id}
^{n-1}\otimes \delta }(S^{n})$ for $n\geq 1$ (see Lemma \ref{leftU(A)}
and Remark \ref{U(A)}) and $U^{0}(S)=\mathbb{C}$.
We first show that the cochain complex that defined the
dual cohomology can be ``restricted to $U^{n}(S)$'' in some cases.

\medskip
\begin{lemma}
Let $(R,\delta)$ be a (not necessarily saturated) Hopf $C^*$-algebra.
Suppose that $(X,\beta ,\gamma )$ is a $R$-bicomodule such that $
\gamma =1_{R}\otimes \mathrm{id}_{X}$. Then $\partial
_{n}(CB(X;U^{n}(R)))\in CB(X;U^{n+1}(R))$.
\end{lemma}
\noindent \textbf{Proof:} Let $T\in \mathrm{CB}(X;U^{n}(R))$, $s\in R$ and $
x\in X$. Then
\[
((\mathrm{id}^{n}\otimes \delta )(T\otimes \mathrm{id})\beta
(x))(1\otimes s)=(T\otimes \mathrm{id}\otimes \mathrm{id})(\beta
\otimes \mathrm{id})(\beta (x)\cdot s)\in (T\otimes \mathrm{id}
)(M_{R}(X\otimes R))\otimes R\subseteq M(R^{n+1})\otimes R
\]
(note that $(T\otimes \mathrm{id}
)\circ \beta \otimes \mathrm{id}$ is a $R$-bimodule map by Lemma \ref
{1.12}(a)) and
\[
((\mathrm{id}^{n}\otimes \delta )(\mathrm{id}\otimes T)\gamma (x))(1\otimes
s)=1\otimes (\mathrm{id}^{n-1}\otimes \delta )T(x)(1\otimes s)\in
M(R^{n+1})\otimes R
\]
(by the definition of $U^n(R)$). Moreover,
\[
((\mathrm{id}^{n}\otimes \delta )(\mathrm{id}^{n-1}\otimes \delta
)T(x))(1\otimes s)=(\mathrm{id}^{n-1}\otimes \delta \otimes \mathrm{id})((
\mathrm{id}^{n-1}\otimes \delta )T(x)(1\otimes s))\in M(R^{n+1})\otimes R
\]
and for $2\leq k\leq n$,
\[
((\mathrm{id}^{n}\otimes \delta )(\mathrm{id}^{n-k}\otimes \delta \otimes
\mathrm{id}^{k-1})T(x))(1\otimes s)=(\mathrm{id}^{n-k}\otimes \delta \otimes
\mathrm{id}^{k})((\mathrm{id}^{n-1}\otimes \delta )T(x)(1\otimes s))\in
M(R^{n+1})\otimes R.
\]
These show that $((\mathrm{id}^{n}\otimes \delta )\partial
_{n}(T)(x))(1\otimes s)\in M(R^{n+1})\otimes R$. Similarly, we also have $
(1\otimes s)((\mathrm{id}^{n}\otimes \delta )\partial _{n}(T)(x))\in
M(R^{n+1})\otimes R$. Thus $\partial _{n}(T)(x)\in U^{n+1}(R)$ as required.

\bigskip

The above lemma says that for any right $R$-bicomodule $(X,\beta )$,
if we take $\gamma =1_{R}\otimes \mathrm{id}_{X}$, then $(\mathrm{CB}
(X;U^{n}(R)),\partial _{n})$ is a cochain subcomplex of $(\mathrm{CB}
(X;M(R^{n})),\partial _{n})$. The cohomology $H_{R,d}^{n}(R;X)$ defined by
this complex is called the \emph{restricted left trivial dual cohomology. }
It turns out that the vanishing of $H_{R,d}^{1}(R;X)$ is related to the
existence of a left invariant mean on $U^{1}(R)$. The idea of
the necessity of part (a) in the following proposition comes from \cite[p.43]
{Pat}.

\medskip
\begin{proposition}
\label{exist-im2}
(a) Suppose that $(R,\delta)$ is a counital
(not necessarily saturated) Hopf $C^{\ast }$-algebra such that
$R$ has property (S) (in particular, if $R$ is a nuclear
$C^{\ast }$-algebra).
Then $H_{R,d}^{1}(R;X)=(0)$ for any right $R$-comodule $X$ if and
only if there exists a left invariant mean $\Phi $ on $U^{1}(R)$
(see Definition \ref{imean}).

\smallskip
\noindent (b) If $(S,\delta)$ is a saturated Hopf $C^{\ast }$-algebra
such that $S$ is unital,
then $H_{R,d}^{n}(S;X)=(0)$ for all $n\geq 1$ and for all right $S$-comodule $X$.
\end{proposition}
\noindent \textbf{Proof:} (a) Let $X=U^{1}(R)/\mathbb{C}\cdot 1$ with the
canonical quotient map $q$ from $U^{1}(R)$ to $X$. As $(q\otimes \mathrm{id}
)\delta (1)=0$ (where $q\otimes\id$ is the map from $M_R(U^1(R)\otimes
R)$ to $M_R(X\otimes R)$ given by Lemma \ref{1.12}(a)),
the coproduct $\delta $ induces a right coaction $\beta $
from $X$ to $M_{R}(X\otimes R)$ such that
$\beta \circ q=$ $(q\otimes \mathrm{id})\circ \delta $
(by Lemma \ref{A30}(c)).
$X$ can now be regarded as a $R$-bicomodule if we take $\gamma =1\otimes
\mathrm{id}_{X}$. Consider $T=\mathrm{id}_{U^{1}(R)}-\epsilon\cdot 1$ (where $\epsilon$ is
the counit of $R$). Since $T(1)=0$, it induces a map $\bar{T}\in \mathrm{
CB}(X;U^{1}(R))$ such that $\bar{T}\circ q=T$ (Lemma \ref{A1}).
Now for any $z\in U^{1}(R)$,
\begin{eqnarray*}
\partial _{1}(\bar{T})(q(z)) &=&(T\otimes \mathrm{id})(\delta (z))-\delta
(T(z))+1\otimes T(z) \\
&=&\delta (z)-1\otimes z-\delta (z)+\epsilon(z)\cdot 1\otimes 1+1\otimes
z-\epsilon(z)\cdot 1\otimes 1\quad =\quad 0.
\end{eqnarray*}
\noindent Hence by the hypothesis, there exists $\phi \in X^{\ast }$ such that $\bar{T}
=\partial _{0}(\phi )$, i.e. $(\phi \circ q\otimes \mathrm{id})\circ \delta
-\phi \circ q\cdot 1=\mathrm{id}_{U^{1}(R)}-\epsilon\cdot 1$. 
Now, let $\Phi =\epsilon-\phi \circ q$.
It is clear that $\Phi (1)=\epsilon(1)=1$ and ($\Phi \otimes 
\mathrm{id})\circ \delta = \mathrm{id}_{U^{1}(R)}-(\phi \circ q\otimes 
\mathrm{id})\circ \delta =\epsilon\cdot 1-\phi \circ q\cdot 1$. 
By a similar argument as \cite[2.2]{Pat} (see also \cite[2.1]{Ruan3}), the 
rescaling of either the positive part or the negative part of $\Phi$ in the Jordan decomposition is a left invariant mean. 
Conversely, suppose that there exists a left
invariant mean $\Phi $ on $U^{1}(R)$. Let $(X,\beta )$ be any right $R$-comodule. For any $T\in \mathrm{CB}(X;U^{1}(R))$ such that $\partial
_{1}(T)=0$, take $f=\Phi \circ T\in X^{\ast }$.
Then for any $g\in R^{\ast }$ and $x\in X$, we have by Lemma \ref{1.12}(c),
\begin{eqnarray*}
0 &=&\Phi ((\mathrm{id}\otimes g)((T\otimes \mathrm{id})\beta
(x)))-\Phi ((\mathrm{id}\otimes g)\delta (T(x)))+\Phi (1)g(T(x)) \\
&=&g((f\otimes \mathrm{id})\beta (x))-\Phi (T(x))g(1)+g(T(x)).
\end{eqnarray*}
\noindent Hence $T(x)=(f\otimes \mathrm{id})\beta (x)-f(x)\cdot 1=\partial _{0}(f)$.

\smallskip
\noindent (b) It was shown in \cite{Wor} (see also \cite{vD0})
that $S$ has a left Haar state $
\phi $ which is in fact a left invariant mean on $U^{1}(S)=S$. Now for any $
T\in \mathrm{Ker}(\partial _{n})$, if $F=(\phi \otimes \mathrm{id}
^{n-1})\circ T$,
\[
0=(F\otimes \mathrm{id})\circ \beta +\sum_{k=1}^{n-1}(-1)^{k}(\mathrm{
id}^{n-k-1}\otimes \delta \otimes \mathrm{id}^{k-1})\circ F+(-1)^{n}1\otimes
F+(-1)^{n+1}T
\]

\noindent which completes the proof.

\bigskip

This proposition, together with Remarks \ref{U(A)}(a) and
\ref{remim}(b), gives the following corollary.
We recall that a discrete
semi-group $\Lambda $ is said to be \emph{left amenable }if there exists a
left invariant mean on $l^{\infty }(\Lambda )$ (see Remark \ref{remim}(b)).
Note that $U(\Lambda )=U_{l}(\Lambda )=l^{\infty }(\Lambda )$ (see Remark
\ref{U(A)}(a)) and right coactions of $c_{0}(\Lambda )$ on an operator space $Z$
are bounded homomorphisms from $\Lambda $ to $\mathrm{CB}(Z;Z)$
(by a similar argument as in Example \ref{B1}(c)).

\medskip
\begin{corollary}
(a) Suppose that $M$ is $a$ locally compact semi-group with identity. Then $
U_{l}(M)$ (see Remark \ref{U(A)}(a)) has a left invariant mean if and only
if $H_{R,d}^{1}(C_{0}(M);Y)=(0)$ for any right $C_{0}(M)$-comodule $Y$.

\smallskip
\noindent (b) If $\Lambda $ is a discrete semi-group with identity, then $
\Lambda $ is left amenable if and only if for any bounded representation $
\pi $ of $\Lambda $ in $\mathrm{CB}(Z;Z)$ and any $T\in \mathrm{CB}
(Z,l^{\infty }(\Lambda ))$ with $T(z)(r\cdot s)=T(\pi(s)z)(r)+
T(z)(s)$ ($r,s\in \Lambda $; $z\in Z$), there exists an element
$f\in Z^{\ast }$ such that $T(z)(r)=f(\pi (r)z)-f(z)$.
\end{corollary}

\medskip
Propositions \ref{D8} and \ref{exist-im2} can also be used to prove the following
interesting theorem.

\medskip
\begin{theorem}
\label{D9} Let G be a locally compact group.

\smallskip
\noindent (a) $G$ is amenable if and only if $H_{d}^{1}(C_{0}(G);X)=(0)$ for
any $C_{0}(G)$-bicomodule $X$.

\smallskip
\noindent (b) $G$ is amenable if and only if $H_{d}^{1}(C_{r}^{\ast
}(G);X)=(0)$ for any $C_{r}^{\ast }(G)$-bicomodule $X$ or equivalently,
$H_{d}^{1}(C_{r}^{\ast}(G);X)=(0)$ for any 2-sided non-degenerate
$C_{r}^{\ast }(G)$-bicomodule $X$.
\end{theorem}
\noindent \textbf{Proof: }(a) If $G$ is amenable, then it is clear that $
OH^{1}(L^{1}(G);X^{\ast })=(0)$ (since the Banach algebra cohomology $
H^{1}(L^{1}(G);X^{\ast })$ vanishes) and thus $L^{1}(G)$ is operator
amenable. Now by putting $S=C_{0}(G)$ and $A=L^{1}(G)$ into Proposition \ref
{D8} (see \cite[\S 5]{Ng1}), we see that $H_{d}^{1}(C_{0}(G);X)=(0)$. The
converse follows directly from Proposition \ref{exist-im2}(a)
(note that for any right $S$-comodule $X$, if we take $\gamma = 1\otimes
\id_X$, then $H^1_{R,d}(S;X)\subseteq H^1_d(S;X)$) as well as
Remarks \ref{U(A)}(a) (i.e. $U(G)\subseteq U^1(C_0(G))$)
and \ref{remim}(b).

\smallskip
\noindent (b) Suppose that $G$ is amenable. Then by
\cite[3.6]{Ruan2}, it is easily seen that $H_{d}^{1}(C_{r}^{\ast}(G);X)$
is zero if we put $S=C_{r}^{\ast }(G)$ and $A=A(G)$ into
Proposition \ref{D8} (see e.g. \cite[2.5]{CR} or \cite[\S 5]
{Ng1}).
On the other hand, if $H_{d}^{1}(C_{r}^{\ast
}(G);X)=(0)$ for any 2-sided non-degenerate
$C_{r}^{\ast }(G)$-bicomodule $X$, then Proposition \ref
{B12}(c) tells us that $C_{r}^{\ast }(G)$ has a counit and hence $G$ is
amenable.

\medskip
\begin{remark}
\label{4.7}(a) If the Fourier algebra $A$ of a saturated Hopf
$C^*$-algebra $S$ (recall that the {\em Fourier algebra} is the
intersection of all non-zero $S$-invariant (closed) ideals of
$S^*$; see \cite[\S 5]{Ng1}) separates points of $S$, then
by the proof of Proposition \ref{D8}, we have the inclusions:
$H_{d}^{1}(S;X)\subseteq
H_{d}^{1}(A^{\ast };X)\cong OH^{1}(A;X^{\ast })\subseteq H^{1}
(A;X^{\ast })$.
In fact, these inclusions, together with the results in \cite{Ruan2} and
\cite{John}, give one of the implications of both parts (a) and (b) of
the above Theorem.
The other implications tell us that even the vanishing of
$H_{d}^{1}(C_{0}(G);X)$ or $H_{d}^{1}(C_{r}^{\ast }(G);X)$
is strong enough to characterise the amenability of $G$.

\smallskip
\noindent (b) The amenability of $G$ is also equivalent to the vanishing of
the first dual cohomology of $vN(G)$ (which is exactly \cite[3.6]{Ruan2}
by Proposition \ref{C15}). The same is true for $L^{\infty }(G)$
(by the inclusions in part (a) as well as the results in \cite{John}).
\end{remark}

\medskip
\begin{corollary}
$\label{cor-sc}$Let $G$ be a locally compact group.

\smallskip
\noindent (a) If $G$ is amenable, then there exists a codiagonal on $
\hat{U}(G\times G)=\{g\in C_{b}(G\times G):r\mapsto r\ast _{1}g$ and $r
\mapsto r\ast _{2}g$ are both continuous maps$\}$ (where $r\ast
_{1}g(s,t)=g(r^{-1}s,t)$ and $r\ast _{2}g(s,t)=\Delta (r^{-1})g(s,tr^{-1})$).

\smallskip
\noindent (b) If there exists a codiagonal $F$ on C$_{b}(G\times G)$ such 
that $F\circ \delta = \epsilon$ on $C_b(G)$,
then $G$ is amenable.
\end{corollary}
\noindent
{\bf Proof:}
(a) This follows from Proposition \ref{exist-sc}(a) and
Theorem \ref{D9}(a) (together with a similar identification as in Remark
\ref{U(A)}(a)).

\smallskip\noindent
(b) For any right $R$-comodule $(X,\beta)$, if we consider $\gamma$ to
be the 2-sided non-degenerate left coaction $1\otimes \id_X$, then
Proposition \ref{exist-sc}(b) implies that $(0) = H^1_d(C_0(G);X)
\supseteq H^1_{R,d}(C_0(G); X)$.
Hence by Proposition \ref{exist-im2}(a), there exists a left invariant
mean on $U^1(C_0(G))\supseteq U(G)$ and $G$ is amenable.

\bigskip

Next, we will consider the vanishing of the one-sided dual cohomology and
relate it to the amenability of Hopf $C^{\ast }$-algebras. All the
remaining results in this section are obvious and hence no proof will be
given.

\medskip

First of all, we have the following one-sided version of Theorem \ref{D9}
which is a direct consequence of Proposition \ref{B12}.

\medskip
\begin{corollary}
\label{Ba} Let $G$ be a locally compact group.

\smallskip
\noindent (a) H$_{d,r}^{1}(C_{0}(G);Y)=(0)$ for any right $C_{0}(G)$-comodule $Y$.

\smallskip
\noindent (b) $G$ is amenable if and only if $H_{d,r}^{n}(C_{r}^{\ast
}(G);Y)=(0)$ for any $n\in \mathbb{N}$ and any right $C_{r}^{\ast }(G)$-comodule $Y$.
It is the case if and only if $H_{d,r}^{1}(C_{r}^{\ast
}(G);Y)=(0)$ for any 2-sided non-degenerate right $C_{r}^{\ast }(G)$-comodule $Y$
\end{corollary}

\medskip
For the general case, we need to consider the cohomology of Hopf von
Neumann algebras instead. We first recall from \cite[2.1]{Ng2} that a
Hopf $C^{\ast } $-algebra $(S,\delta)$ is said to be \emph{left
H}$_{1}$\emph{-coamenable} if the
{\em left Fourier algebra} $A_{S}^{l}$ (i.e. the intersection of all non-zero $S$-invariant left (closed) ideals of $S^{\ast }$; see \cite[\S 5]{Ng1}) has a
bounded left approximate
identity. By Theorem \ref{C6}, we have a characterisation of the amenability
of $S$ in terms of the cohomology of the Hopf von Neumann algebra $
(A_{S}^{l})^{\ast }$.

\medskip
\begin{corollary}
\label{D6} Suppose that the left Fourier algebra $A$ of $(S,\delta)$ is non-zero.
Then $(S,\delta)$ is left $H_{1}$-coamenable if and only if $H_{d,r}^{1}
(A^{\ast };X)=(0)$ for any right $A^{\ast }$-comodule $X$.
\end{corollary}

\medskip
Let $(T,V,S)$ be a Kac-Fourier duality in the sense of \cite[5.13]{Ng1}.
Then by \cite[3.15]{Ng2}, the left (or right) $H_{1}$-coamenability of $T$ will
automatically imply the 2-sided $H_{1}$-coamenability (i.e. the Fourier
algebra of $S$ as defined in Remark \ref{4.7}(a) has a bounded 2-sided
approximate identity). In this case, we will simply call $(T,V,S)$ \emph{
amenable} (see \cite[3.16]{Ng2}). We have the following characterisation of
this amenability in terms of cohomology.

\medskip
\begin{proposition}
\label{D3} Let $(T,V,S)$ be a Kac-Fourier duality and $(\mu ,\nu )$ be
any $V$-covariant representation on a Hilbert space $H$. Let $B=\mu ^{\ast }(
\mathcal{L}(H)_{\ast })$ (note that $B$ will then be both the Fourier algebra and the
left Fourier algebra of $S$ by \cite[5.9]{Ng1}). Then $(T,V,S)$ is amenable
if and only if $H_{d,r}^{1}(B^{\ast };X)=(0)$
for any right $B^{\ast }$-comodule $X$.
\end{proposition}

\medskip
Note that the cohomology considered in this proposition can be regarded as
an one-sided version of the operator cohomology (using Proposition \ref{C15}).

\medskip

Suppose that $V\in \mathcal{L}(H\otimes H)$ and $W\in \mathcal{L}(K\otimes K)$
are regular multiplicative unitaries such that $S_{V}\cong \hat{S}_{W}$
(see \cite[1.5]{BS}), in particular, if $V$ comes from a Kac system $(H,V,U)$
(see \cite[\S 6]{BS}) and $W = \Sigma (1\otimes U)V(1\otimes U)\Sigma$
(where $\Sigma \in \mathcal{L}(H\otimes H)$ is defined by 
$\Sigma(\xi\otimes\eta) = \eta\otimes\xi$).
Then $(\hat{S}_{V},V,S_{V})$ is a Kac Fourier
duality (see \cite{Ng1}). In this case, $V$ is amenable in the sense of Baaj
and Skandalis (see \cite[A13(c)]{BS}) if and only if $(\hat{S}_{V},V,S_{V})$
is amenable in the above sense.
Moreover, if $\hat{A}_{V}=L_{V}^{\ast }(\mathcal{L}(H)_{\ast })$
(where $L_{V}$ is the default representation of $S_{V}$ on
$\mathcal{L}(H)$), then $\hat{A}_{V}^{\ast }\cong S_{V}^{\prime \prime }$
(the weak*-closure of $S_{V}$ in $
\mathcal{L}(H)$).
Hence we can express amenability of Kac systems (see
\cite[3.3\&6.2]{BS}) and in particular Kac algebras (see \cite{ES2}) in
terms of cohomology.

\medskip
\begin{corollary}
\label{D5} Suppose that $V$ is a regular irreducible multiplicative unitary.
$V$ is amenable if and only if $H_{d,r}^{1}(S_{V}^{\prime \prime };$ X)=(0)
for any right $S_{V}^{\prime \prime }$-comodule $X$.
\end{corollary}

\medskip
In particular, we have a Hopf von Neumann algebra analogue of Corollary \ref
{Ba}(b) which can be regarded as an one-sided version of \cite[3.6]{Ruan2}
(in the light of Proposition \ref{C15}).

\bigskip

\bigskip

\appendix

\section{Extensions of coactions}

\bigskip
\noindent
{\bf Notation:}
{\em In this appendix, we do not assume the Hopf $C^*$-algebra
$(R,\delta)$ to be saturated.
Moreover, as usual, all the right and left coactions on $C^*$-algebras
in this appendix are $*$-homomorphisms.}

\bigskip
The main objective of this appendix is to answer the
following natural and interesting question: ``\textit{If} $\beta $
\textit{is a coaction of} $R$ \textit{on a}
$C^{\ast }$-algebra $A$\textit{, what is the biggest unital
closed subalgebra of }$M(A)$
\textit{on which }$\beta $
\textit{can be extended?}'' This extension is important for some
arguments in this paper. We know that in general, it is impossible to
extend $\beta $ to a
coaction on the whole of $M(A)$ (consider e.g. $A=C_{0}(G)=R$ and $
\beta $ is the coproduct on $C_{0}(G)$). Inspired by the
definition of $U(G)$, we have the following lemma.

\medskip
\begin{lemma}
\label{leftU(A)} Suppose that $R$ has property (S) in the sense of \cite{Was}
or \cite{AB} (in particular, if it is a nuclear $C^{\ast }$-algebra; see
\cite[10]{Was}). Let $\beta $ be a (right)
coaction of $R$ on a $C^*$-algebra $A$.
Then $U_{\beta }(A)=\{m\in M(A):\beta (m)\in M_{R}(M(A)\otimes R)\}$
is a unital $C^{\ast }$-subalgebra of $M(A)$ and $\beta $ extends to a
coaction (again denoted by $\beta $) of $R$ on U$_{\beta }(A)$.
\end{lemma}
\noindent
\textbf{Proof:} It is obvious that $U_{\beta }(A)$ is a
unital $C^{\ast }$-subalgebra of $M(A)$.
We need to show the inclusion: $\beta
(U_{\beta }(A))\subseteq M_{R}(U_{\beta }(A)\otimes R)$.
In fact, since $\beta (U_{\beta }(A))(1\otimes R)\subseteq
M(A)\otimes R$ and $R$
has property (S), it suffices to prove that $(\mathrm{id}\otimes f)
\beta (m)\in U_{\beta }(A)$ for any $f\in R^{\ast }$ and $
m\in U_{\beta }(A)$.
Note that there exist $f^{\prime }\in R^{\ast }$ and
$t\in R$ such that $f=t\cdot f^{\prime }$.
Thus for any $s\in R$, we have
$\beta \lbrack (\mathrm{id}\otimes f)\beta (m)](1\otimes s)=
(\mathrm{id}\otimes \mathrm{id}\otimes f^{\prime})
[(\id\otimes \delta)\beta (m)(1\otimes s\otimes t)]\in
M(A)\otimes R$
(recall that $(\id\otimes \delta)\beta(m)\in M_{R\otimes R}
(M(A)\otimes R \otimes R)$ by Lemma \ref{1.12}(b)).
Finally, the coaction identity of the extension
follows from that of $\beta $.

\bigskip

It is clear that $U_{\beta }(A)$ is the biggest closed
subalgebra of $M(A)$ on which $\beta $ can be extended. Similarly, we
can define $U_{\gamma }(A)$ for a left coaction $\gamma $.
Moreover, we have the following two-sided version of Lemma \ref{leftU(A)}.

\medskip
\begin{lemma}
\label{2sU(A)} Suppose that $R$ has property (S) and there exists
a (right) coaction $\beta $ and a left coaction $\gamma $ of $R$ on a
$C^{\ast }$-algebra $A$ such that $(\gamma\otimes \id)\beta = (\id\otimes
\beta)\gamma$.
Then $U_{\beta ,\gamma }(A)=\{m\in M(A):\beta (m)\in
M_{R}(M(A)\otimes R)$; $\gamma (m)\in M_{R}(R\otimes M(A))\}$ is a
unital $C^{\ast }$-subalgebra of $M(A)$ and $\beta $
(respectively, $\gamma $) extends
to a (right) coaction (respectively, left coaction), again denoted by $\beta $
(respectively, $\gamma $), on $U_{\beta ,\gamma }(A)$.
Moreover, $U_{\beta ,\gamma }(A)$ is the biggest unital $C^*$-subalgebra
of $M(A)$ for which both $\beta$ and $\gamma$ can be extended.
\end{lemma}
\noindent \textbf{Proof:} Note that $U_{\beta ,\gamma
}(A)=U_{\beta }(A)\cap U_{\gamma }(A)$ and so is a unital
$C^{\ast }$-subalgebra of $M(A)$. We first show that $\beta $ extends to a
coaction on $U_{\beta ,\gamma }(A)$. By Lemma \ref{leftU(A)} and its
proof, we need only to show that $(\mathrm{id}\otimes f)\beta (m)\in U_{\gamma
}(A)$ for any $m\in U_{\beta ,\gamma }(A)$
and $f\in R^{\ast }$.
Indeed, for any $s\in R$, we have $\gamma ((\mathrm{id}\otimes f)\beta
(m))(s\otimes 1)=(\mathrm{id}\otimes \mathrm{id}\otimes f)(\mathrm{id}
\otimes \beta )(\gamma (m)(s\otimes 1))\in R\otimes M(A)$ and
similarly $(s\otimes 1)\gamma ((\mathrm{id}\otimes f)\beta (m))\in
R\otimes M(A)$. The proof for the extension of $\gamma $ is the same.

\medskip
\begin{remark}
\label{U(A)}
If $A=R$ and $\beta =\delta =\gamma $, we denote $U_{\beta
,\gamma }(A)$ by $U(R)$. Moreover, we denote $U_{\mathrm{id}\otimes \delta
,\delta \otimes \mathrm{id}}$($R\otimes R)$ by $\hat{U}(R\otimes R)$ and $U_{
\mathrm{id}^{n-1}\otimes \delta }(R^{n})$ by $U^{n}(R)$ ($n\geq 1$).

\smallskip\noindent
(a) Suppose that $M$ be a locally compact semi-group. Then $
(C_{0}(M),\delta _{M})$ (where $\delta _{M}(f)(r,s)=f(rs)$) is a (possibly
non-saturated) Hopf $C^{\ast }$-algebra. Let $U(M)=\{g\in C_{b}(M):$ $
r\mapsto r\cdot g$ and $r\mapsto g\cdot r$ are both norm continuous$\}$ and $
U_{l}(M)=\{f\in C_{b}(M):$ $r\mapsto r\cdot f$ is a norm continuous map$\}$
(where $r\cdot f(s)=f(sr)$ and $f\cdot r(s)=f(rs))$. Then $
U_{l}(M)=U^{1}(C_{0}(M))$ and $U(M)=U(C_{0}(M))$.
In fact, for any $f\in C_{b}(M)$, using a similar argument as in Example
\ref{B1}(c), $\delta (f)\in M_{C_{0}(M)}(C_{b}(M)\otimes
C_{0}(M))=C_{b}(M;C_{b}(M))$ if and only if the map that sends
$r\in M$ to r$\cdot f\in C_{b}(M)$ is continuous.
Consequently, for a locally compact group $G$, the space $U(C_{0}(G))$
coincides with $U(G)$ (note that the inverse and the modular
function are not included in the definition of $U(R)$ but it doesn't matter).

\smallskip
\noindent (b) Let $R=C_{r}^{\ast }(G)$ and $u_{t}$ be the element in $
M(C_{r}^{\ast }(G))$ corresponding to $t\in G$. Then it is clear that $
u_{t}\in U(R)$.
\end{remark}

\medskip
\begin{definition}
\label{imean} Suppose that $\beta $ is a (right) coaction of $R$ on a
$C^{\ast}$-algebra $A$.

\smallskip
\noindent (a) A closed subspace $X$ of M(A) is said to be \emph{weakly }$
\beta $\emph{-invariant} if $(\id\otimes
f)\beta (X)\subseteq X$ for any $f\in R^{\ast }$.

\smallskip
\noindent (b) Let $X$ be a weakly $\beta $-invariant subspace of $M(A)$
that contains 1$_{A}$. Then $\Phi \in X^{\ast }_+$ is said to be a $
\beta $-\emph{invariant mean} if $\Phi (1_{A})=1$ and
$\Phi ((id\otimes f)\beta (x))=\Phi (x)f(1_{R})$
for any $f\in R^{\ast }$ and $x\in X$.
In the case when $A=R$ and $\beta =\delta $, we call
such $\Phi $ a \emph{left invariant mean }on $X$.
\end{definition}

\medskip

We can define similarly weakly $\gamma $-invariant
subspaces and $\gamma $-invariant mean for a left coaction $\gamma $.

\medskip
\begin{remark}
\label{remim}(a) It is
clear that if $\beta$ induces a right coaction on a subspace $X$
of $M(A)$ (i.e.
$\beta (X)\subseteq M_{R}(X\otimes R)$), then $X$ is automatically weakly $
\beta $-invariant.
If in addition $1_A\in X$, then $\Phi \in X^{\ast }_+$
is a left invariant mean if $\Phi (1_{A})=1$ and $(\Phi \otimes \mathrm{id}
)\circ \beta =\Phi \cdot 1_{A}$ on $X$. This applies, in particular, to
both $U^{1}(R)$ and $U(R)$ if $R$ has property (S) (by Lemmas \ref{leftU(A)}
and \ref{2sU(A)}).

\smallskip
\noindent (b) If $M$ is a locally compact semi-group, then $(C_{0}(M),\delta
_{M})$ is a (not necessarily saturated) Hopf $C^{\ast }$-algebra. Let $X$ be
a \emph{left} \emph{invariant} subspace of $
C_{b}(M)$ in the sense that $r\cdot g\in X$
for any $g\in X$ and $r\in M$ where $r\cdot g(t)=g(tr)$.
It is clear that $\delta _{M}$ induces a right coaction on $X$
and so $X$ is weakly $\delta _{M}$-invariant (by part (a)). If $X$ contains $
1$, then $\Phi \in X^{\ast }_+$ is a left invariant mean in the sense
of Definition \ref{imean} if and only if $\Phi (1)=1$ and $\Phi (r\cdot
f)=\Phi (f)$ for any $f\in X$. Hence for a locally compact group $G$, the
left invariant mean on $U(G)=U(C_0(G))$ as defined above coincides with the usual
definition.
\end{remark}

\medskip

\noindent
Department of Pure Mathematics,
The Queen's University of Belfast,
Belfast BT7 1NN,
Northern Ireland,
United Kingdom

\smallskip\noindent
E-mail address: c.k.ng@qub.ac.uk


\bigskip

\bigskip

\begin{thebibliography}{99}
\bibitem{AB}  R.J. Archbold and C.J.K. Batty, $C^{\ast }$-tensor norms and
slice maps, J. Lond. Math. Soc.(2), 22 (1980), 127-138.

\bibitem{BS0}  S. Baaj and G. Skandalis, $C^{\ast }$-alg\`{e}bres de Hopf et
th\'{e}orie de Kasparov \'{e}quivariante, K-theory, 2 (1989), 638-721.

\bibitem{BS}  S. Baaj and G. Skandalis, Unitaires multiplicatifs et dualit
\'{e} pour les produits crois\'{e}s de $C^{\ast }$-alg\`{e}bres, Ann.
scient. \'{E}c. Norm. Sup., $4^{e}$ s\'{e}rie, t. 26 (1993), 425-488

\bibitem{Blech}  D.P. Blecher, The standard dual of an operator space,
Pacific J. Math., 153 (1992), 15-30.

\bibitem{BP}  D.P. Blecher and V.I. Paulsen, Tensor products of operator
spaces, J. Funct. Anal., 99 (1991), 262-292.

\bibitem{CES}  E. Christensen, E.G. Effros and A. Sinclair, Completely
bounded multiplinear maps and $C^{\ast }$-algebraic cohomologoy, Invent.
Math., 90 (1987), 279-296.

\bibitem{CM}  A. Connes and H. Moscovici, Cyclic Cohomology, Hopf Algebras
and the modular theory, preprint: math.QA/9905013.

\bibitem{CR}  Jean De Canni\`{e}re and Ronny Rousseau, The Fourier Algebra
as an Order Ideal of the Fourier-Stieltjes Algebra, Math. Z., 186 (1984),
501-507.

\bibitem{Doi}  Y. Doi, Homological coalgebra, J. Math. Soc. Japan, 33
(1981), 31-50.

\bibitem{EKR}  E. Effros, J. Kraus and Z.-J. Ruan, On two quantized tensor
products, \textit{Operator algebras, mathematical physics, and
low-dimensional topology (Istanbul, 1991)}, Res. Notes Math., 5 (1993),
125-145.

\bibitem{ES}  M. Enock and J.-M. Schwartz, Algebres de Kac Moyennables,
Pacific J. Math., 125 (1986), 363-379.

\bibitem{ES2}  M. Enock and J.-M. Schwartz, \textit{Kac algebras and duality
of locally compact groups}, 1992, Springer-Verlag.

\bibitem{Evens}  L. Evens, \textit{The Cohomology of Groups}, Oxford
Mathematical Monogrophs, 1991, Clarendon Press.

\bibitem{John}  B.E. Johnson, \emph{Cohomology in Banach algebras}, Mem.
Amer. Math. Soc. 127, 1972.

\bibitem{John2}  B.E. Johnson, Non-amenability of the Fourier algebra of a
compact group, J. Lond. Math. Soc., 50 (1994), 361-374.

\bibitem{Muhly}  P.S. Muhly, A finite-dimensional introduction to Operator
algebras, \emph{Operation algebras and Applications (Samos, 1996)}, SE: NATO
Adv. Sci. Inst. Ser. C Math. Phys. Sci. 495, 313-354. \noindent

\bibitem{Ng0}  C.K. Ng, Discrete coaction on $C^{\ast }$-algebras, J.
Austral. Math. Soc. (Series A), 60 (1996), 118-127.

\bibitem{Ng}   C.K. Ng, Coactions and crossed products of Hopf $C^*$-algebras,
Proc. London Math. Soc.(3) 72(1996), 638-656.

\bibitem{Ng1}  C.K. Ng, Duality of Hopf $C^{\ast }$-algebras, preprint.

\bibitem{Ng2}  C.K. Ng, Amenability of Hopf $C^{\ast }$-algebras,
Proceedings of the 17th conference on Operator Theory, to appear.

\bibitem{Ng3}  C.K. Ng, Coactions on operator spaces and exactness,
unpublished work.

\bibitem{Ng4}  C.K. Ng, Profinite Quantum groups, in preparation. 

\bibitem{Pal}  T.W. Palmer, Banach algebras and the general theory of 
$\ast $-algebras Vol. I: Algebras and Banach algebras, Encyclopedia of Math. 1994,
Camb. Univ. Press.

\bibitem{Pat}  A. T. Paterson, \textit{Amenability}, Mathematical Surveys
and Monographs no. 29, 1991, Amer. Math. Soc.

\bibitem{Ped}  G. K. Pedersen, $C^{\ast }$\emph{-algebras and their
automorphism groups}, 1979, Academic Press.

\bibitem{Pier}  J.-P. Pier, \emph{Amenable Banach algebras}, Pit. Res. Note.
Math. Series 172, 1988, Longman.

\bibitem{Ruan}  Z. J. Ruan, On the predual of dual algebras, J. Oper.
Theory, 27 (1992), 179-192.

\bibitem{Ruan2}  Z. J. Ruan, The operator amenability of $A(G)$, Amer. J.
Math., 117 (1995), 1449-1474.

\bibitem{Ruan3} Z. J. Ruan, Amenability of Hopf von Neumann algebras and 
Kac algebras, J. Funct. Anal. 139 (1996), 466-499.

\bibitem{Serre}  J.-P. Serre, \textit{Galois Cohomology}, 1997,
Springer-Verlag, Berlin.

\bibitem{Shn}  S. Shnider, Bialgebra deformations, C. R. Acad. Sci. Paris,
S\'{e}rie I, t. 312 (1991), 7-12.

\bibitem{SS}  A.M. Sinclair and R.R. Smith, \emph{Hochschild cohomology of
von Neumann algebras}, Lond. Math. Soc. LNM 203, 1995, Cambridge University
Press.

\bibitem{Take}  M. Takesaki, \textit{Theory of Operator algebras I}, 1979,
Springer-Verlag.

\bibitem{vD0} A. van Daele, The Haar measure on a Compact Quantum
Group, Proc. Amer. Math. Soc. 123 (1995), 3125-3128.

\bibitem{VW}  A. Van Daele and S. Z. Wang, Universal quantum groups,
Internat. J. Math., 7 (1996), no. 2, 255-264.

\bibitem{Wang}  S.Z. Wang, Free Products of Compact Quantum Groups,
Comm. Math. Phys. 167 (1995), 671-692.

\bibitem{Was}  S. Wassermann, A pathology in the ideal space of $L(H)\otimes
L(H)$, Indiana Univ. Math. J., 27 (1978), 1011-1020.

\bibitem{Wilson}  J. Wilson, \textit{Profinite groups}, LMS Monographs New
Series 19, 1998, Clarendon Press.

\bibitem{Wor}  S. L. Woronowicz, Compact matrix pseudogroups, Comm. Math.
Phys., 111(1987), 613-665.

\end{thebibliography}
\end{document}